\documentclass[leqno]{amsart}
\usepackage{amsmath}
\usepackage{amsthm}
\usepackage{amssymb}
\usepackage{latexsym}
\usepackage{enumerate}

\newtheorem{thm}{Theorem}[section]
\newtheorem{prop}{Proposition}[section]
\newtheorem{lem}{Lemma}[section]
\newtheorem{cor}{Corollary}[section]
\newtheorem{rem}{Remark}[section]
\newtheorem{defi}{Definition}[section]

\usepackage{enumerate}

\allowdisplaybreaks

\title{Limiting absorption principle on manifolds having ends with various measure growth rate limits}

\author{Hironori Kumura${}^*$}

\email{smhkumu@ipc.shizuoka.ac.jp} 
\address{Department of Mathematics, Shizuoka University, 
Shizuoka 422-8529, Japan}

\thanks{${}^*$\ The author thanks Kenichi Ito for helpful discussions. }

\date{revised form, Nov., 2011}

\begin{document}

\maketitle

\begin{abstract}
The purpose of this paper is to study the property of the resolvent of the Laplace-Beltrami operator on a noncompact complete Riemannian manifold with various ends each of which has a different limit of the growth rate of the Riemannian measure at infinity, in particular, focusing on the limiting absorption principle. 
As a result, we will obtain the absolute continuity of the Laplace-Beltrami operator. 
\end{abstract}
\section{Introduction}

The Laplace-Beltrami operator $\Delta_g$ on a noncompact complete Riemannian manifold $(M,g)$ is essentially self-adjoit on $C^{\infty}_0(M)$, and its self-adjoit extension to $L^2(M, v_g)$ has been studied by several authors from various points of view. 
The purpose of this paper is to prove the limiting absorption principle and the absolute continuity of the Laplace-Beltrami operator on a Riemannian manifold with a combination of a plurality of ends, each of which has a {\it different limit of the growth rate of the Riemannian measure} at infinity. 
A Riemannian manifold with a plurality of ends has been studied by several authors (see Mazzeo-Melrose \cite{M-M}, Mazzeo \cite{Ma}, Perry \cite{P}, Bouclet \cite{B} and so on), but they assumed that the curvatures of {\it all} ends converged to a {\it common constant} at infinity; such an assumption does not produce any difference between ``one'' and ``more-than-one'' with respect to the number of ends. 
This paper studies a Riemannian manifold with a combination of different geometries for each end. 

It is important to note that the spectral structure of the Laplace-Beltrami operator determines the asymptotic behavior of a solution of the wave equation or of a time-dependent Schr\"odinger equation as time tends to infinity. 
For example, if $u$ is a vector of the absolutely continuous subspace of the Laplace-Beltrami operator, then the ``wave function'' $e^{i\Delta t}u$ decays locally as $t \to \pm \infty$; on the other hand, if $u$ is an element of the pure point subspace of the Laplace-Beltrami operator, the wave function $e^{i\Delta t}u$ remains localized for all time $t$. 

To state our results, we shall first introduce some terminology and notations. 
Let $(M,g)$ be an $n$-dimensional connected noncompact complete Riemannian manifold and $V$ be an (possibly unbounded) open subset of $M$. 
We shall say that $E:=M \backslash V$ is an {\it end with radial coordinates} if and only if $V$ has a compact connected $C^{\infty}$-boundary $\partial V$ such that the outward normal exponential map $\exp_{\partial V}^{\perp} : N^{+}(\partial V) \to M \backslash  V$ induces a diffeomorphism, where $N^{+}(\partial V) := \{w \in T(\partial V) \mid w~{\rm is~outward~normal~to}~\partial V \}$. 
We assume that there exists a relatively compact open subset $U$ of $M$ such that $M \backslash U$ consists of a finite number of ends $E_1, E_2, \cdots , E_m$ with radial coordinates;
\begin{align*}
  M \backslash U = E_1 \cup E_2 \cup \cdots \cup E_m~~({\rm disjoint~union}),
\end{align*}
where $m \ge 1$ is an integer. 
We denote $r(x):= {\rm dist}_g(U,x)$ for $x \in M \backslash U$. 
For convenience, we shall introduce the following terminologies: we shall say that, an end $E$ with radial coordinates satisfies the condition ``MC$\displaystyle \big( \frac{a}{r}; \frac{b}{r}; \delta \big)$'', if there exist constants $a>0$, $b>0$, and $\delta \in (0,1)$ such that
\begin{align}
 \begin{cases}
  \ \displaystyle \nabla dr 
  \ge \Big\{ \frac{a}{r} + O (r^{-1-\delta}) \Big\}( g - dr \otimes dr ) 
  \quad & {\rm on}~~E, \vspace{2mm} \\
  \ \displaystyle \Delta_g r = \frac{b}{r} + O (r^{-1-\delta}) 
  \quad & {\rm on}~~E; 
 \end{cases}
\end{align}
also, we shall say that an end $E$ satisfies the condition ``MC$(\alpha; \beta; \delta)$'', if there exist constant $\alpha > 0$, $\beta > 0$, and $\delta \in (0,1)$ such that
\begin{align}
 \begin{cases}
  \ \displaystyle \nabla dr 
  \ge \big\{ \alpha + O (r^{-1-\delta}) \big\}( g - dr \otimes dr ) 
  \quad & {\rm on}~~E, \vspace{2mm} \\
  \ \displaystyle \Delta_g r = \beta + O (r^{-1-\delta}) 
  \quad & {\rm on}~~E; 
 \end{cases}
\end{align}
``MC'' stands for ``mean curvature''; note that $\Delta_g r$ is the mean curvatures of each level-hypersurface of the function $r$, and that $\Delta_g r$ expresses the growth rate of the Riemannian measure $v_g$ on the level-hypersurfaces of $r$. 
We assume that 
\begin{align}
 & E_j \ \ {\rm satisfies} \ \ 
   {\rm MC} \Big( \frac{a_j}{r}; \frac{b_j}{r}; \delta_j \Big) 
   \quad {\rm for}~~ 1 \le j \le m_0 ;\\
 & E_j \ \ {\rm satisfies} \ \ 
   {\rm MC}(\alpha_j ; \beta_j ; \delta_j) 
   \qquad \,{\rm for}~~ m_0 + 1 \le j \le m ,
\end{align}
where $a_1, a_2, \cdots , a_{m_0}$, $b_1, b_2, \cdots , b_{m_0}$, $\alpha_{m_0+1}, \alpha_{m_0+2}, \cdots , \alpha_{m}$ are positive constant; 
\begin{align*}
  0 = \beta_1 = \beta_2 = \cdots = \beta_{m_0} < \beta_{m_0+1} 
  \le \beta_{m_0+2} \le \cdots \le \beta_m
\end{align*}
are real constants; $\delta_1, \delta_2, \cdots, \delta_m$ are constant in $(0,1)$; $0$ ``$\le$'' $m_0 \le m$ is an integer. 
\begin{rem}
{\rm If $m_0=0$, we mean that every end satisfies $(4)$, and that there is no end satisfying $(3)$, and hence, $\{ a_1, a_2, \cdots ,a_{m_0} \} = \{ b_1, b_2, \cdots , b_{m_0} \} = \{ \delta_1, \delta_2, \cdots , \linebreak \delta_{m_0} \} = \emptyset$.} 
\end{rem}
In order to state our theorems, we need to introduce the weighted $L^2$-spaces; let $v_g$ be the Riemannian measure of $(M,g)$. 
For any $s \in \mathbb{R}$, let $L^2_s(M,?v_g)$ denote the Hilbert space of all complex-valued measurable functions $f$ such that $\left|(1+r)^s f\right|$ is square integrable on $M$; the inner product and norm will be denoted as follows:
\begin{align*}
 & (f, g)_{L^2_s(M, v_g)} = \int_M \left\{ 1 + r(x) \right\}^{2s} f(x)\,
    \overline{g(x)} \,dv_g(x);\\
 & \| f \|_{L^2_s(M, v_g)} = \sqrt{(f, f)_{L^2_s(M, v_g)}},
\end{align*}
where recall that $r={\rm dist}_g(U,*)$. 
Let $R(z)$ denote the resolvent $( - \Delta_g - z )^{-1}$ of $-\Delta_g$ for $z\in \rho(-\Delta_g)$, where $\rho(-\Delta_g)$ stands for the resolvent set of $-\Delta_g$. 
For $t, t' \in \mathbb{R}$, let $\mathbb{B}(t, t')$ denote the space of all bounded linear operators $T: L^2_{t}(M, v_g) \to L^2_{t'}(M, v_g)$ with ${\rm Dom}\,(T) = L^2_{t}(M, v_g)$; 
$\mathbb{B}(t, t')$ is a Banach space with the operator norm $\| * \|_{\mathbb{B}(t, t')}$. 
We set 
\begin{align*}
 \delta := \min \left\{ \delta_j \mid 1\le j\le m \right\}
\end{align*} 
and denote
\begin{align*}
  I 
& := \left( \frac{(\beta_1)^2}{4} , \, \infty \right) - 
  \left\{ \frac{(\beta_j)^2}{4} \ \Big| \ j = 2, 3, \cdots, m \right\};\\
  {\it \Pi}_{+} 
& := \left\{ z = x + iy \in \mathbb{C} \ \Big| \ x > \frac{(\beta_1)^2}{4}, \, y \ge 0 \right\} 
  - \left\{ \frac{(\beta_j)^2}{4} \ \Big| \  j= 2,3,\cdots ,m \right\};\\
  {\it \Pi}_{-} 
& := \left\{ z = x + iy \in \mathbb{C} \ \Big| \  x > \frac{(\beta_1)^2}{4}, \, y \le 0 \right\} 
  - \left\{ \frac{(\beta_j)^2}{4} \ \Big| \  j= 2,3,\cdots ,m \right\}.
\end{align*}
Now, we shall state our main result. 
\begin{thm}[{\bf Principle of limiting absorption}]
Let $(M, g)$ be an $n$-dimensional connected complete Riemannian manifold. 
Assume that there exists a relatively compact open subset $U$ of $M$ such that $M \backslash U$ consists of a disjoint union of finitely many ends $E_1, E_2, \cdots , E_m$ with radial coordinates, where $m \ge 1$ is an integer. 
Assume that the conditions $(3)$ and $(4)$ holds, where $a_1, a_2, \cdots , a_{m_0}$, $b_1, b_2, \cdots , b_{m_0}$, $\alpha_{m_0+1}, \alpha_{m_0+2}, \cdots , \alpha_{m}$ are positive constant\,{\rm ;} 
\begin{align*}
  0 = \beta_1 = \beta_2 = \cdots = \beta_{m_0} < \beta_{m_0+1} 
  \le \beta_{m_0+2} \le \cdots \le \beta_m
\end{align*}
are real constants\,{\rm ;} $\delta_1, \delta_2, \cdots, \delta_m$ are constant in $(0,1)$\,{\rm ;} $0$ ``$\le$'' $m_0 \le m$ is an integer. 
Note that, if $m_0 = 0$, we follows the convention stated in {\rm Remark} $1.1$. 
Let $s$ and $s'$ be real numbers satisfying
\begin{align*}
  0 < s' < s < \min \Big\{ a_{\min}, \frac{1}{2} \Big\}
  ; \quad s' + s \le \delta = \min \{ \delta_1, \cdots, \delta_m \},
\end{align*}
where 
\begin{align*}
 a_{\min} := 
 \begin{cases}
  \ \min \{ a_1, \cdots, a_{m_0} \} \qquad & {\rm if}~~\beta_1=0;\\
  \ \infty \qquad & {\rm if}~~\beta_1>0. 
 \end{cases}
\end{align*}
Then, in the Banach space $\mathbb{B} \big( \frac{1}{2} + s, - \frac{1}{2} - s' \big) $, we have the limit 
\begin{align}
  R (\lambda \pm i \,0)
  := \lim_{\varepsilon \downarrow 0} R(\lambda \pm i \varepsilon)
  \qquad {\rm for}~~\lambda \in I \qquad ({\rm double~sign~in~same~order}).
\end{align}
Moreover, the convergence $(5)$ is uniform on any compact subset of $I$, and hence, $R(z)$ is continuous on ${\it \Pi}_{+}$ and ${\it \Pi}_{-}$ with respect to the operator norm $\| * \|_{ \mathbb{B} \big( \frac{1}{2} + s, - \frac{1}{2} - s' \big) }$, by considering $R(\lambda + i\,0)$ on ${\it \Pi}_{+} \cap (0,\infty)$, and $R(\lambda - i\,0)$ on ${\it \Pi}_{-} \cap (0,\infty)$, respectively.
\end{thm}
In order to state the next theorem, we shall recall some terminology from the spectral theory: let $H$ be a self-adjoint operator on a Hilbert space $(X, (,)_X)$, and $E(\Lambda)$ $(\Lambda \in \mathcal{B})$ denote the spectral measure of $H$ on $X$, where $\mathcal{B}$ stands for the set of all Borel sets of $\mathbb{R}$; ``spectral measure'' is also called ``spectral decomposition'' or ``spectral projection''. 
Any $u \in X$ defines the measure $m_u$ on $\mathbb{R}$, by $m_u(\Lambda) := ( E(\Lambda)u, u)_{X}$ for $\Lambda \in \mathcal{B}$. 
Let $|*|$ denote the Lebesgue measure on $\mathbb{R}$. 
The property of $m_u$ classifies vectors in $X$ as follows: set 
\begin{align*}
  X_{\rm pp} 
& := {\rm the~closure~of~the~linear~hull~of~eigenvectors~of}~H ;\\
  X_{\rm c} 
& := \left\{ u \in X \mid m_u(\{a\})= 0~
  {\rm for~any}~a \in \mathbb{R} \right\}; \\
  X_{\rm ac} 
& := \bigl\{ u \in X_{\rm c} \mid 
  ~m_u~{\rm is~absolutely~continuous~with~respect~to}~|*| \bigr\}; \\
  X_{\rm sc} 
& := \bigl\{ u \in X_{\rm c} \mid 
  ~m_u~{\rm is~singular~continuous~with~respect~to}~|*| \bigr\};
\end{align*}
then, $X$ is decomposed into the direct sum of three closed linear subspaces:
\begin{align}
  X = X_{\rm pp} \oplus X_{\rm ac} \oplus X_{\rm sc},
\end{align}
where $X_{\rm pp}$, $X_{\rm ac}$, and $X_{\rm sc}$ are orthogonal to each other, and reduce $H$; thus, corresponding to the decomposition $(6)$, $H$ is decomposed into the direct sum of three self-adjoint operators:
\begin{align*}
  H = H_{\rm pp} \oplus H_{\rm ac} \oplus H_{\rm sc}.
\end{align*}
$X_{\rm pp}$, $X_{\rm c}$, $X_{\rm ac}$, and $X_{\rm sc}$ are called the {\it pure point subspace} of $H$, {\it continuous subspace} of $H$, {\it absolutely continuous subspace} of $H$, and {\it singular continuous subspace} of $H$, respectively. 
If $H = H_{\rm ac}$ (i.e., $X=X_{\rm ac}$), then $H$ is said to be {\it absolutely continuous}. 
Also, if $H = H_{\rm sc}$ (i.e., $X=X_{\rm sc}$), $H$ is said to be {\it singular continuous}. 
Moreover, for an open interval $J$ in $\mathbb{R}$, if $H|_{E(J)X}$ is absolutely continuous, then $H$ is said to be {\it absolutely continuous on $J$}. 
(For the details mentioned above, see, \cite{Ka} Chapter $10$, or \cite{R-S} VII.$2$). 
\begin{thm}
Assume that $(M, g)$ satisfies the assumptions in Theorem $1.1$. 
Then, $\sigma_{\rm ess}(-\Delta_g) = [ \frac{(\beta_1)^2}{4} , \infty )$, and $-\Delta_g$ is absolutely continuous on $( \frac{(\beta_1)^2}{4} , \infty )$, where $\sigma_{\rm ess} (-\Delta_g) $ stands for the essential spectrum of $-\Delta_g$. 
In particular, $-\Delta_g$ has no singular continuous spectrum. 
\end{thm}
\begin{cor}
Assume that $(M, g)$ satisfies the assumptions in Theorem $1.1$. 
If $m_0\ge 1$, then $-\Delta_g$ is absolutely continuous on $(0,\infty)$ and $0$ is not an eigenvalue of $-\Delta_g$. 
\end{cor}
Note that the decay order $(2)$ is fairly sharp; indeed, there exists a rotationally symmetric manifold $( \mathbb{R}^n, g := dr^2 + f(r)^2 g_{S^{n-1}(1)})$ which satisfies $\Delta_g r = (n-1) (1+O(r^{-1}))$ as $r \to \infty$; $\sigma_{\rm ess}(-\Delta_g)=[\frac{(n-1)^2}{4}, \infty)$; $\sigma_{\rm pp}(-\Delta_g) \cap ( \frac{(n-1)^2}{4}, \infty ) = \{ \frac{(n-1)^2}{4} + 1 \}$ (see \cite{K3}). 

Now, we shall recall some earlier works and compare them to the theorems above. First, recall that Xavier proved the following:
\begin{thm}[Xavier \cite{X}]
Let $(M,g)$ be an Hadamard manifold and $r$ be the distance function to some fixed point of $M$. 
Assume that the function $f=(r^2+1)^{1/2}$ satisfies the following conditions {\rm (i)} and {\rm (ii):}
{\rm 
\begin{enumerate}[(i)]
 \item  $\Delta_g f \le C_1,$
 \item  $(\Delta_g) ^2 f \le C_2 f^{-3}$,
\end{enumerate}
}\noindent 
where $C_1$ and $C_2$ are positive constants. 
Then, $-\Delta_g$ is absolutely continuous on $(\alpha,\infty)$, where $4\alpha = 6C_1 + C_2$. 
\end{thm}
The nature of Theorem $1.3$ seems to be more analytic than geometric; note that the term $(\Delta_g)^2 f$ in (ii) contains the derivatives of the curvature tensor. 
Theorem $1.1$ and $1.2$ seem to be more geometric than Xavier's in the sense that Theorem $1.1$ and $1.2$ do not require any estimates of derivatives of the curvature tensor. 

After that, Donnelly proved the following by using the Mourre theory:
\begin{thm}[Donnelly \cite{D3}]
Assume that a complete Riemannian manifold $(M,g)$ admits a proper $C^2$-exhaustion function $b$ satisfying the following\,{\rm :}
{\rm 
\begin{enumerate}[(i)]
 \item  $c_1 r \le b \le c_2 r$\, for some positive constants 
$c_1$ and $c_2$$;$
 \item  $\left||\nabla b| - 1 \right| \le c \,b^{-\varepsilon}$$;$
 \item  $\left| \nabla db - \displaystyle b^{-1}(g-db\otimes db) \right|
	 \le c\,b^{-1-\varepsilon}$$;$
 \item  $\left|(b^2)_{kks}\right| + \left| (b^2)_{skk} \right| \le 
        c\,b^{-\varepsilon}$,
\end{enumerate}
}
\noindent 
where $c$ and $\varepsilon$ are positive constants{\rm ;} $r$ denotes the distance function to a fixed point of $M${\rm ;} indices stand for the components of the covariant differential{\rm ;} the repeated indices expresses contraction. 
Then $-\Delta_g $ is absolutely continuous. 
\end{thm}

Since the function $b$ is not a distance function, Theorem $1.4$ seems to be general in that sense. 
But, the author has a feeling that the condition (iv) of Theorem $1.4$ is somehow technical. 

In this paper, we will modify the arguments used in a classical method for the Schr\"odinger operator in Euclidean space; see, Eidus \cite{E1} \cite{E2}, Ikebe and Saito \cite{I-S}, and Mochizuki and Uchiyama \cite{M-U} and so on. 
We will dare to use the distasnce function explisitely. 

The Mourre theory (\cite{Mo}, \cite{A-B-G}) is a powerful pool for studying of the Schr\"odinger operator on the Euclidean space, and its nature is quite abstract. Note that Froese and Hislop \cite{F-H} studied spectral properties of second-order elliptic operators on noncompact manifolds by using the Mourre theory from ``analytic'' points of view (they required the estimates of some derivatives of the curvature tensor); Gol\'{e}nia and Moroianu \cite{G-M} proved the limiting absorption principle under a bounded condition on the second derivative of the metric on conformally cusp manifolds by modifying the Mourre theory; see also Froese-Hislop-Perry \cite{F-H-P} for a hyperbolic manifold, and Guillop\'{e} \cite{G} and so on. 
\vspace{2mm}


{\bf Acknowledgements.} 
The author would like to express his gratitude to Minoru Murata, Hiroshi Isozaki, Kenichi Ito, and Sylvain Gol\'{e}nia for useful discussions. 

\section{Unitarily equivalent operator $L$ and radiation condition}

In this section, we shall define the unitarily equivalent operator $L$, and introduce the radiation condition for $L$. 

First, we shall list the notation used in the sequel: 
\begin{align*}
  & E_j(s,t) := \{ x\in E_j \mid s < r(x) < t \};\\
  & E_j(s,\infty) := \{ x\in E_j \mid s < r(x)\};\\
  & S_j(t) := \{ x \in E_j \mid r(x) = t \};\\
  & E(s,t) := \{ x\in M \mid s<r(x)<t\} 
    = \displaystyle \cup_{j=1}^{m}E_j(s,t);\\
  & E(s,\infty) := \{ x \in M \mid s < r(x) \}
    = \displaystyle \cup_{j=1}^{m} E_j(s,\infty);\\
  & S(t) := \{ x \in M \mid r(x) = t \}
    = \displaystyle \cup_{j=1}^{m}S_j(t);\\
  & U(R) := U \cup \{ x \in M \mid 0 \le r(x) < R \}, 
\end{align*}
where $R > 0$; $0 \le s < t$. \\

We take a real-valued $C^{\infty}$ function $w$ on $M$ so that 
\begin{align*}
  w(x) = \frac{\beta_j}{2} \,r(x) \quad {\rm for} \quad x \in E_j(1,\infty) 
  \quad (1\le j\le m),
\end{align*}
where recall $r(x)={\rm dist}_g(U,x)$ for $x \in M\backslash U$. 
We shall introduce a new measure $\mu$ on $M$ and the operator $L$ as follows:
\begin{align}
  \mu := & e^{-2w}v_g; \quad L := e^w \circ \Delta_g \circ e^{-w}.
\end{align} 
Then, a direct computation shows that 
\begin{align}
 Lf = \Delta_g f - 2\langle \nabla w, \nabla f \rangle -Vf \,; \quad 
 V := \Delta_g w - |\nabla w|^2.
\end{align}
Since the multiplication operator $e^w:L^2(M, v_g) \ni h \mapsto e^{w}h \in L^2(M,\mu)$ is unitary, $L$ with $\mbox{Dom}(L) = \{ u \in L^2(M,\mu) \mid Lu \in L^2(M,\mu) \}$ is a nonnegative self-adjoint operator on $L^2(M,\mu)$. 
Note that assumptions $(3)$ and $(4)$, together with $(1)$ and $(2)$, imply that $\Delta_g r \to \beta_j$ as $r \to \infty$ on $E_j$, and hence, 
\begin{align}
  V(x) = \frac{\beta_j}{2} \, \Delta_g r(x) - \Big( \frac{\beta_j}{2} \Big)^2 
  \to  \frac{(\beta_j)^2}{4} \quad {\rm as}~~x \in E_j~
  {\rm and}~r(x)\to \infty.
\end{align}
In particular, $V$ is bounded on $M$. 

Let $A$ denote the induced Riemannian measures on each level hypersurface $S(t)$ for $t \ge 0$, and set 
\begin{align*}
  A_w := e^{-2w} A.
\end{align*}
For $\Omega \subset M$ and $s \in \mathbb{R}$, let $L^2_s(\Omega ,\mu)$ denote the Hilbert space of all complex-valued measurable functions $f$ such that $|(1+r)^sf|$ is square integrable over $\Omega$. 

Now, we shall consider the equation 
\begin{align}
  - Lu - zu = f
\end{align}
for some suitably chosen $z \in \mathbb{C}$ and $f\in L^2(M,\mu )$. 
If $z\in \rho(-L)$, this equation has a unique solution $u\in L^2(M,\mu )$ for $f\in L^2(M,\mu )$, where $\rho(-L)$ is the resolvent set of $-L$. 
In order to extend this uniqueness theorem for $z\in \sigma_{{\rm c}}(-L)$, we have to consider the operator $-L$ in a wider class and introduce the boundary condition at infinity for our manifolds, where $\sigma_c(-L)$ stands for the continuous spectrum of $-L$; we choose smooth functions $p_{+} : M \times {\it \Pi}_{+} \to \mathbb{C}$ and $p_{-} : M \times {\it \Pi}_{-} \to \mathbb{C}$ so that, for each $1\le j\le m$,  
\begin{align}
  p_{\pm}(x,z) = \mp i \sqrt{z - \frac{(\beta_j)^2}{4} }~r(x) + \alpha_j \log r(x) 
  \quad {\rm for}~~(x,z)\in E_j(r_0,\infty)\times \Pi_{\pm},
\end{align}
where $r_0 \ge 1$ is a constant; $\alpha_j$ is the function defined by
\begin{equation}
 \alpha_j(x) :=
 \begin{cases}
  \ \displaystyle \frac{b_j}{2} \quad & \text{if} \ 
      x \in E_{j}(r_0,\infty) \quad ( 1\le j \le m_0 ), \vspace{1mm}\\
  \ \ 0 \quad & \text{if} \ 
      x \in E_{j}(r_0,\infty) \quad ( m_{0}+1 \le j \le m ).
 \end{cases}
\end{equation}
Note that the square root in $(11)$ is the principal value, that is, the analytic extension of $\sqrt{x}$ for $x > 0$. 
Note also that we are following the usual ``double sign in same order" convention'' in $(11)$. 
This convention will be used in the sequel of this paper. 
We shall consider the following condition:
\begin{defi}[{\bf radiation condition}]
{\rm We shall say that a solution $u$ of the equation $(10)$ with $z \in {\it \Pi}_{\pm}$ satisfies the {\it radiation condition} if there exists constants $s'$ and $s$ such that 
\begin{align}
 & 0<s'\le s <1;~~s'+ s \le 1;\\
 & u \in L^2_{- \frac{1}{2} - s' }(M, \mu),
   ~~\partial_r u + (\partial_r p_{\pm}) u \in L^2_{-\frac{1}{2}+s}(M,\mu).
\end{align}
A solution $u$ of $(10)$ satisfying the radiation condition will be called an {\it outgoing solution} or {\it incoming solution}, if $z \in {\it \Pi}_{+}$ or $z \in {\it \Pi}_{-}$, respectively.}
\end{defi}

\section{Energy integral}

This section will be devoted to proofs of Proposition $3.1$ and $3.2$ below, which express energy integrals of a solution of $(10)$. 

First, we extend the Riemannian metric $g = \langle *, * \rangle$ to the complex bilinear form for complex tangent vectors: $ \langle u_1 + i v_1 , u_2 + i v_2 \rangle = \langle u_1 , u_2 \rangle - \langle v_1 , v_2 \rangle + i \bigl\{ \langle u_1 , v_2 \rangle + \langle v_1 , u_2 \rangle \bigr\}$ for $u_1,u_2,v_1,v_2 \in T_xM$~~$(x\in M)$; we also denote $ | u + i v |^2 = \langle u, u \rangle + \langle v, v \rangle $ for $u, v \in T_xM~~(x\in M)$.

For the sake of convenience of readers, we mention two lemmas below, which will immediately follow from the standard Green's formula and divergence theorem, respectively: 
\begin{lem}[Green's formula]
Let $\Omega $ be a relatively compact open subset of $M$ with $C^{\infty}$-boundary $\partial \Omega $, and $u$ and $v$ be $C^{\infty}$-functions on $M$. 
Then, we have
\begin{align*}
  \int_{\Omega} (Lu)v \,d\mu 
= \int_{\partial \Omega }\langle \nabla u, \overrightarrow{n} \rangle v \,dA_w
  - \int_{\Omega} \langle \nabla u, \nabla v \rangle \,d\mu
  - \int_{\Omega} Vuv \,d\mu,
\end{align*}
where $\overrightarrow{n}$ stands for the outward unit normal vector field along $\partial \Omega$. 
\end{lem}
\begin{lem}[divergence theorem]
Let $\Omega $ be a relatively compact open subset of $M$ with $C^{\infty}$-boundary $\partial \Omega $, and $X$ be a $C^{\infty}$-vector field on $M$. 
Then, we have
\begin{align*}
  \int_{\Omega} ( {\rm div}\,X) \,d\mu 
= \int_{\partial \Omega } \langle X, \overrightarrow{n} \rangle \,dA_w
  + 2 \int_{\Omega} \langle X, \nabla w \rangle \,d\mu.
\end{align*}
\end{lem}
\begin{prop}
Let $\varphi(r)$ be a nonnegative function of $r\ge 0$ and $u\in H^2_{loc}(M)$ be a solution of the equation $(10)$. 
Let $\Omega $ be a relatively compact open subset of $M$ with $C^{\infty}$-boundary $\partial \Omega $. 
Then, we have 
\begin{align*}
 & - \int_{\partial \Omega} \varphi 
   \langle \nabla u, \overrightarrow{n} \rangle \overline{u} \,dA_w
   - z \int_{\Omega } \varphi |u|^2 \,d\mu \\
=& \int_{\Omega } \left\{ \varphi \left( f \overline{u} - |\nabla u|^2 
   - V |u|^2 \right) - \varphi ' 
   \left( \partial_r u \right) \overline{u} \right\} \,d\mu.
\end{align*}
\end{prop}
\begin{proof}
Multiplying the equation $(10)$ by $\varphi \overline{u}$ and integrating it over $\Omega$ with respect to the measure $\mu$, we obtain the desired equation by Lemma $3.1$. 
\end{proof}
The following proposition will play an important role in our arguments:
\begin{prop}
Let $\varphi=\varphi(r)$ be a real-valued function of $r\in [0,\infty)$ satisfying $\varphi(r) \ge 0$ for $r>0$, and $u \in H^2_{{\rm loc}}(M)$ be a solution of $(10)$ satisfying the radiation condition. 
Then, for any $R > r_0$, we have 
\begin{align*}
 & \int_{S(R)} \varphi |{\rm Im}\,\partial_r p_{\pm}| |u|^2  \,dA_w 
   + |{\rm Im}\,z| \int_{U(R)} \varphi |u|^2 \,d\mu \\
\le 
 & \int_{S(R)} \varphi \bigl| {\rm Im}\,\bigl( \overline{u} 
   (\partial_r + \partial_r p_{\pm})u \bigr) \bigr| \,dA_w 
   + \int_{U(R)} \Bigl\{ \varphi |{\rm Im}\,(f\overline{u})| 
   + |\varphi'{\rm Im}\,( \overline{u} \, \partial_r u )| \Big\} \,d\mu,
\end{align*}
where we set $( \partial_r + \partial_r p_{\pm})u := \partial_r u + (\partial_r p_{\pm}) u$ for simplicity. 
Here, note that 
\begin{align*}
  {\rm Im}\,\partial_r p_{\pm} 
  = \mp {\rm Re}\,\sqrt{z - \frac{(\beta_j)^2}{4} } 
  \quad {\rm on~each}~E_j(r_0,\infty).
\end{align*}
\end{prop}
\begin{proof}
Applying Proposition $3.1$ for $U(R)$, we obtain
\begin{align*}
 & \int_{S(R)} \varphi (\partial_r p_{\pm}) |u|^2 \,dA_w 
   - z \int_{U(R)} \varphi |u|^2 \,d\mu \\
=& \int_{S(R)} \varphi \overline{u} (\partial_r + \partial_r p_{\pm})u \,dA_w 
   + \int_{U(R)} \Bigl\{ \varphi \big( f \overline{u} - |\nabla u|^2-V|u|^2 
   \big) -  \varphi' (\partial_r u)\overline{u} \Big\} \,d\mu.
\end{align*}
Taking the imaginary part of both sides of this equation, we get
\begin{align}
 & \int_{S(R)} \bigl( {\rm Im}\,(\partial_r p_{\pm}) \bigr) \varphi |u|^2 \,dA_w   - ({\rm Im}\,z) \int_{U(R)} \varphi |u|^2 \,d\mu \\
=& \int_{S(R)} \varphi \,{\rm Im}\, \bigl( \overline{u} (\partial_r 
   + \partial_r p_{\pm})u \bigr) \,dA_w + \int_{U(R)} \big\{ \varphi \,{\rm Im}
   \,(f \overline{u}) - \varphi' \,{\rm Im}\, \bigl( (\partial_r u)\overline{u}    \bigr) \bigr\} \,d\mu. \notag 
\end{align}
Note that, for a general $z' \in \mathbb{C}$, 
\begin{align*}
 {\rm Re}\,\sqrt{z'} \quad 
 \begin{cases}
  > 0  \quad & {\rm if}~~ {\rm Im}\,z' > 0, \\
  > 0  \quad & {\rm if}~~ z' \in (0,\infty), \\
  = 0  \quad & {\rm if}~~ z' \in (-\infty,0], 
 \end{cases}
\end{align*}
where recall that we take the principal value as our square root. 
Hence, by $(11)$ and $(12)$, we see that, 
\begin{align}
& {\rm if}~z \in {\it \Pi}_{+}, \hspace{2mm} 
  {\rm Im}\,z \ge 0 \hspace{2mm} {\rm and} \hspace{2mm} 
  {\rm Im}\,\partial_r p_{+} 
  = -{\rm Re}\,\sqrt{z - \frac{(\beta_j)^2}{4} } \le 0~\,{\rm on}~E_j(r_0,\infty); \\
& {\rm if}~z \in {\it \Pi}_{-}, \hspace{2mm} 
  {\rm Im}\,z \le 0 \hspace{2mm} {\rm and} \hspace{2mm} 
  {\rm Im}\,\partial_r p_{-} 
  = {\rm Re}\,\sqrt{z - \frac{(\beta_j)^2}{4} } \ge 0~\,{\rm on}~E_j(r_0,\infty). 
\end{align}
Thus, signs of ${\rm Im}\,\partial_r p_{+}$ and ${\rm Im}\,z$ are different; moreover, $\varphi |u|^2 \ge 0$. 
Hence, Proposition $3.2$ follows from $(15)$, $(16)$, and $(17)$. 
\end{proof}
In the sequel, we will simply write ${\rm Im}\, (\partial_r u)\overline{u} :={\rm Im}\, \bigl( (\partial_r u)\overline{u} \bigr)$ and so on. 
\section{A priori estimate of $|\nabla u|$}

We shall introduce an operator $L_{\rm loc}$ by ${\rm Dom}\,(L_{\rm loc}) = H^2_{\rm loc}(M)$ and $L_{\rm loc}u := \Delta u - 2\langle \nabla w, \nabla u \rangle - Vu$ for $u \in H^2_{\rm loc}(M)$.
Then, the following holds: 
\begin{lem}[local a priori estimate]
Let $\Omega$ be a domain of $M$ and $\varphi$ be a real-valued function of $r\ge 0$. 
Assume that 
${\rm supp} \,\varphi$ is compact\,{\rm ;} $\varphi\big|_{\partial \Omega} = 0$\,{\rm ;} $|\varphi| \le 1$. 
Moreover, assume that there exist a constant $s \in \mathbb{R}$ and a function $u$ such that $u \in H^2_{{\rm loc}}(\Omega) \cap L^2_s( \Omega, \mu)$ and $L_{\rm loc}u \in L^2_s( \Omega, \mu)$. 
Then, for any $\varepsilon\in (0,1)$, we obtain
\begin{align*}
 & ( 1 - \varepsilon ) \int_{\Omega} \varphi^2 (1+r)^{2s} |\nabla u|^2 
   \,d\mu\\
\le 
 & \,\frac{\varepsilon }{2} \int_{\Omega} \varphi^2 (1+r)^{2s} 
   |L_{{\rm loc}}u|^2 \,d\mu 
   + \widehat{c}_0 \int_{{\rm supp}~\varphi} (1+r)^{2s} |u|^2 \,d\mu,
\end{align*}
where 
\begin{align*}
 \widehat{c}_0 := \frac{1}{2\varepsilon} + \max_M |V| + \frac{1}{\varepsilon} 
 \max_M \left|\varphi' + s (1+r)^{-1} \varphi \right|^2.
\end{align*}
\end{lem}
\begin{proof}
Since $\varphi ^2 ( 1 + r )^{2s} \nabla u = \nabla \left\{ \varphi ^2 ( 1 + r )^{2s} u \right\} - 2u \varphi ( 1 + r )^{2s} \{ \varphi ' + s (1+r)^{-1} \varphi \} \nabla r$, Lemma $3.1$ implies that 
\begin{align*}
 & \int_{\Omega} \varphi^2 (1+r)^{2s} |\nabla u|^2 \,d\mu \\ 
=& \int_{\Omega} \langle \nabla \left\{ \varphi ^2 ( 1 + r )^{2s} u \right\},      \nabla \,\overline{u} \rangle \,d\mu 
   - 2 \int_{\Omega} \varphi ( 1 + r )^{2s} 
   \{ \varphi ' + s (1+r)^{-1} \varphi \} u \,\partial_r \overline{u} \,d\mu \\
=& - \int_{\Omega} ( L_{\rm loc} \overline{u} ) \varphi ^2 ( 1 + r )^{2s} u 
   \,d\mu - \int_{\Omega} \varphi ^2 ( 1 + r )^{2s} V |u|^2 \,d\mu \\
 & - 2 \int_{\Omega} \varphi ( 1 + r )^{2s} \{ \varphi ' + s (1+r)^{-1} 
   \varphi \} u \,\partial_r \overline{u} \,d\mu \\
\le 
 & \frac{\varepsilon}{2} \int_{\Omega} \varphi ^2 ( 1 + r )^{2s} 
   |L_{\rm loc} u|^2 \,d\mu 
   + \frac{1}{2\varepsilon} \int_{\Omega} \varphi ^2 ( 1 + r )^{2s} |u|^2 
   \,d\mu \\
 & + \max_M |V| \cdot \int_{\Omega} \varphi ^2 ( 1 + r )^{2s} |u|^2 
   \,d\mu + \varepsilon \int_{\Omega} \varphi ^2 ( 1 + r )^{2s} 
   \left| \partial_r u \right|^2 \,d\mu \\
 & + \frac{1}{\varepsilon} 
   \max_M \left|\varphi ' + s (1+r)^{-1} \varphi \right|^2 
   \cdot \int_{{\rm supp}\, \varphi} ( 1 + r )^{2s} |u|^2 \,d\mu. 
\end{align*}
Now, Lemma $4.1$ follows from this inequality and the assumption $|\varphi| \le 1$. 
\end{proof}
\begin{cor}[global a priori estimate]
Assume that there exit a constant $s\in \mathbb{R}$ and a function $u$ such that $u \in {\rm Dom}(L_{\rm loc}) \cap L^2_s(M,\mu)$ and $L_{\rm loc}u \in L^2_s(M,\mu)$. 
Then, for any $\varepsilon \in (0, \frac{1}{2})$, there exists a constant $\widehat{c}(\varepsilon) > 0$ such that 
\begin{align*}
  \| \nabla u \|^2_{L^2_s(M,\mu)} 
  \le \varepsilon \| L_{{\rm loc}}u \|^2_{L^2_s(M,\mu)}
  + \widehat{c}(\varepsilon) \| u \|^2_{L^2_s(M,\mu)},
\end{align*}
where $\widehat{c}(\varepsilon) := \frac{1}{1-2\varepsilon} \left\{ \frac{1}{4\varepsilon} + \frac{1}{2\varepsilon} ( 1+|s| )^2 + \max_M |V| \right\}$. 
\end{cor}
\begin{proof}
For $t>0$, set
\begin{align*}
 h_t(r) := 
 \begin{cases}
  \ \ \ \ \ 1 \qquad & \mbox{if}\quad r \le t,\\
  - r + t + 1 \qquad & \mbox{if}\quad t \le r\le t+1,\\
  \ \ \ \ \ 0 \qquad & \mbox{if}\quad t+1 \le r
 \end{cases}
\end{align*}
and put $\Omega =M$ and $\varphi(r) = h_t(r)$ in Lemma 4.1. 
Then, 
\begin{align*}
  & (1-2\varepsilon) \int_M h_t^2 (1+r)^{2s} |\nabla u|^2 \,d\mu \\
\le 
  & \varepsilon \int_M h_t^2 (1+r)^{2s} |L_{{\rm loc}}u|^2 \,d\mu 
    + c_t \int_{U(t+1)} (1+r)^{2s} |u|^2 \,d\mu,
\end{align*}
where $c_t := \frac{1}{4\varepsilon} + \max_M |V| + \frac{1}{2\varepsilon} \max_M \left| h_t' + s (1+r)^{-1} h_t \right|^2$. 
Since $|h_t'| \le 1$ and $|s|(1+r)^{-1}|h_t| \le |s|$, we have $c_t \le \frac{1}{4\varepsilon} + \frac{1}{2\varepsilon} ( 1+|s| )^2 
+ \max_M |V|$. 
Letting $t\to \infty$, we get the desired inequality. 
\end{proof}
\section{Estimate of $|\nabla u + u\nabla p_{\pm}|$}
The purpose of this section is to prove Proposition $5.2$ below, which shows a decay estimate of $|\nabla u + u\nabla p_{\pm}|$ of a solution $u$ of $(10)$; decay assumptions $(3)$ and $(4)$ (see also $(1)$ and $(2)$) will be systematically used; in this sense, this is the most important section. 

To prove Proposition $5.2$, we first prove a preparative proposition. 
Let $\eta :M\times {\it \Pi}_{\pm}\to \mathbb{C}$ be a complex-valued $C^{\infty}$-function and consider a function
\begin{align*}
  v(x,z) := e^{\eta (x,z)}u(x),
\end{align*} 
where $u$ is a solution of the equation $(10)$. 
In view of $(8)$ and $(10)$, direct computations show the following:
\begin{align}
 & - \Delta_g v + 2 \langle \nabla \eta + \nabla w, \nabla v \rangle - qv 
   =e^{\eta}f;\\
 & q = q(x,z)
   := z - \Delta_g \eta 
   + \langle 2 \nabla w + \nabla \eta, \nabla \eta \rangle - V.
\end{align}
The following Proposition $5.1$ will serve the estimate of $|\nabla u + u\nabla p_{\pm}|$ on $M$ (see Proposition $5.2$): 
\begin{prop}
Let $\eta : M\times {\it \Pi}_{\pm} \to \mathbb{C}$ and $\psi : M \times {\it \Pi}_{\pm} \to \mathbb{R}$ be $C^{\infty}$-functions and $X$ be a ``real'' $C^{\infty}$-vector field on $M$. 
Let $u$ be a solution of the equation $(10)$ and set $v = e^{\eta }u$. 
Let $\Omega \subset M$ be a relatively compact open subset with $C^{\infty}$-boundary $\partial \Omega$. 
Then, we have 
\begin{align*}
 & \int_{\partial \Omega} \psi \left\{ {\rm Re} \,\langle X, 
   \nabla \overline{v} \rangle \langle \nabla v, \overrightarrow{n} \rangle
   - \frac{1}{2} |\nabla v|^2 \langle X, \overrightarrow{n} \rangle 
   \right\} \,dA_w \\
=& \int_{\Omega} \left( \langle X, \nabla w \rangle \psi 
   - \frac{1}{2} \langle X, \nabla \psi \rangle - \frac{1}{2} 
   \psi \,{\rm div}\,X \right) |\nabla v|^2 \,d\mu \\
 & + \int_{\Omega} \Bigl\{ {\rm Re} \, \langle \nabla \psi 
   + 2 \psi \nabla \eta, \nabla v \rangle \langle X, 
   \nabla \overline{v} \rangle
   + \psi \,{\rm Re} \,\langle \nabla_{\nabla v} X, \nabla \overline{v} \rangle
   \Bigr\} \,d\mu \\
 & - \int_{\Omega} \psi \,{\rm Re} \,( qv + e^{\eta}f )
     \langle X, \nabla \overline{v} \rangle \,d\mu.
\end{align*}
\end{prop}
\begin{proof}
Direct calculations show the following:
\begin{align*}
   - (\Delta_g v) \langle X , \nabla \overline{v} \rangle 
=& - {\rm div}\,\bigl( \langle X, \nabla \overline{v} \rangle \nabla v \bigr) 
   + \langle \nabla_{\nabla v} X, \nabla \overline{v} \rangle 
   + (\nabla d \overline{v})(X,\nabla v) ; \\
   2{\rm Re}\,(\nabla d \overline{v})(X,\nabla v) 
=& \,\,{\rm div}\,\bigl( \langle \nabla v, \nabla \overline{v} \rangle X \bigr)    - \langle \nabla v, \nabla \overline{v} \rangle \, {\rm div}\,X . 
\end{align*}
Combining these equation makes
\begin{align*}
 & - {\rm Re}\,(\Delta_g v) \langle X , \nabla \overline{v} \rangle \\
=& - {\rm Re}\,{\rm div}\,\bigl( \langle X, \nabla \overline{v} \rangle 
   \nabla v \bigr) 
   + {\rm Re}\, \langle \nabla_{\nabla v} X, \nabla \overline{v} \rangle 
   + \frac{1}{2} {\rm div} \, 
   \bigl( | \nabla v |^2 X \bigr) 
   - \frac{1}{2} |\nabla v|^2\, {\rm div}\,X .
\end{align*}
Therefore, multiplying the equation $(18)$ by $\psi \langle X, \nabla \overline{v} \rangle$ and taking its real part, we obtain
\begin{align*}
& - {\rm Re}\, {\rm div} \bigl( \psi \langle X, \nabla \overline{v} \rangle 
  \nabla v \bigr) + {\rm Re}\,\langle X, \nabla \overline{v} \rangle 
  \langle \nabla v, \nabla \psi \rangle 
  + \psi \,{\rm Re} \,\langle \nabla_{\nabla v} X, \nabla \overline{v} 
  \rangle \\
& + \frac{1}{2}\Bigl\{ {\rm div}\, \bigl( \psi |\nabla v|^2 X \bigr) 
  - |\nabla v|^2 \langle X, \nabla \psi \rangle \Bigr\}
  - \frac{1}{2} \psi |\nabla v|^2 {\rm div}\, X \\
& + \psi \,{\rm Re} \Bigl\{ \Bigl( 2\langle \nabla \eta + \nabla w, \nabla v 
  \rangle - qv \Bigr) \langle X, \nabla \overline{v} \rangle \Bigr\} 
  = \psi \,{\rm Re} \,e^{\eta} f \langle X, \nabla \overline{v} \rangle,
\end{align*}
where note that $\psi$ is real-valued. 
Integrating this equation on $\Omega$ with respect to the measure $\mu$ and applying Lemma $3.2$ to the first and fourth term above make
\begin{align*}
 & - \int_{\partial \Omega} \psi \left\{ {\rm Re} \,\langle X, \nabla 
   \overline{v} \rangle \langle \nabla v, \overrightarrow{n} \rangle 
   - \frac{1}{2} |\nabla v|^2 \langle X, \overrightarrow{n} \rangle \right\} 
   \,dA_w \\
 & + 2 \int_{\Omega} \psi \Bigl\{ - {\rm Re}\,\langle X, \nabla \overline{v} 
   \rangle \langle \nabla v, \nabla w \rangle 
   + \frac{1}{2} |\nabla v|^2 \langle X, \nabla w \rangle \Bigr\} \,d\mu \\
 & - \frac{1}{2} \int_{\Omega} |\nabla v|^2 \bigl\{ \langle X, \nabla \psi 
   \rangle + \psi \,{\rm div}\,X \bigr\} \,d\mu 
   + \int_{\Omega} \psi \,{\rm Re} \,\langle \nabla_{\nabla v} X, \nabla 
   \overline{v} \rangle \,d\mu \\
 & + \int_{\Omega} {\rm Re}\,\langle X, \nabla \overline{v} \rangle 
   \Bigl\{ \langle \nabla v, \nabla \psi \rangle 
   + 2 \psi \langle \nabla \eta + \nabla w, \nabla v \rangle - \psi qv 
   - \psi e^{\eta} f \Bigr\} \,d\mu = 0.
\end{align*}
Since the term $2\int_{\Omega}\psi {\rm Re}\,\langle X, \nabla \overline{v} 
   \rangle \langle \nabla v, \nabla w \rangle \,d\mu$ appears twice with different signs on the left hand side of the equation above, we see that Proposition $5.1$ follows from this equation. 
\end{proof}
In the sequel, let $K_{+}$ and $K_{-}$ be any fixed compact subsets in ${\it \Pi}_{+}$ and  ${\it \Pi}_{-}$, respectively. 
Then, we have the following: 
\begin{lem}
If we set $\eta = p_{\pm}$ and if $z\in K_{\pm}$, then the function $q$ defined by $(19)$ has the following asymptotic property on $M$$:$
\begin{align*}
  q = O\left( {r^{-1-\delta}} \right) \quad {\rm uniformly~for}~z\in K_{\pm}~
  {\rm as}~ r \to \infty.
\end{align*}
Here, recall that $\delta = \min\{\delta_j~|~1\le j\le m\}$. 
\end{lem}
\begin{proof}
First, we shall consider the case that $x \in E_j(r_0,\infty)~(1\le j\le m_0)$ and $r(x)\to \infty$. 
Then, $w = 0$; $V = 0$; $\eta = p_{\pm} = \mp i \sqrt{z}~r + \frac{b_j}{2} \log r$, and hence, 
\begin{align}
 q 
 & = z - \Delta_g \eta + \langle 2 \nabla w + \nabla \eta, \nabla \eta \rangle 
   - V  = z - \Delta_g \eta + \langle \nabla \eta, \nabla \eta \rangle \\
 & = \pm i \sqrt{z} \left\{ \Delta_g r - \frac{b_j}{r} \right\}
   - \frac{b_j}{2r}\Delta_g r + \frac{b_j}{2r^2} + \frac{(b_j)^2}{4r^2}. \notag
\end{align}
The condition $(3)$ is $\Delta_g r = \frac{b_j}{r} + O \left( r^{-1-\delta} \right)$ on $E_j$ $(1\le j\le m_0)$, and hence, $(20)$ implies the desired result for ends, $E_1, \cdots, E_{m_0}$. 

Next, we shall consider the case that $x\in E_j(r_0,\infty)~(m_0+1\le j\le m)$ and $r(x)\to \infty$. 
Then, 
$w = \frac{\beta_j}{2} r$; $V = \frac{\beta_j}{2} \Delta_g r - \frac{(\beta_j)^2}{4}$; $\eta = p_{\pm} = \mp i\sqrt{z - \frac{(\beta_j)^2}{4}}\,r$, and hence,
\begin{align}
 q = 
 & z - \Delta_g \eta + \langle 2\nabla w + \nabla \eta, \nabla \eta \rangle 
   - V \\
=& z \pm i \sqrt{ z - \frac{(\beta_j)^2}{4} } \,\Delta_g r 
   + \Big\langle \Big( \beta_j \mp i \sqrt{ z - \frac{(\beta_j)^2}{4} } 
   \,\,\Big) \nabla r, \, 
   \mp i \sqrt{ z - \frac{(\beta_j)^2}{4} }\,\nabla r \Big\rangle \notag \\
 & - \frac{\beta_j}{2} \Delta_g r + \frac{(\beta_j)^2}{4} \notag \\
=& ( \Delta_g r - \beta_j ) \left\{ \pm i \sqrt{ z - \frac{(\beta_j)^2}{4} } 
   - \frac{\beta_j}{2} \right\} \notag.
\end{align}
Since the condition $(4)$ is $\Delta_g r = \beta_j + O\left( r^{-1-\delta} \right)$ on $E_j$ $(m_0 + 1 \le j \le m)$, we see that $(21)$ implies the desired result for ends, $E_{m_0+1}, \cdots, E_{m}$.
\end{proof}
The purpose of this section is to prove the following: 
\begin{prop}
Let $s$, $s'$, and $R$ be positive real numbers, and $z$ be a complex number satisfying 
\begin{align*}
  s < a_{\rm min} = \min \{a_j \mid 1 \le j \le m_0 \},~~
  s + s' \le \delta, \quad 2\le R, \quad z\in K_{\pm}.
\end{align*}
Assume that $f\in L^2_{\frac{1}{2}+s}(M,\mu)$ and that $u$ is a solution of $(10)$ satisfying the radiation condition. 
Then, we have 
\begin{align*}
 & \left\| \nabla u + u \nabla p_{\pm} \right\|^2_{ L^2_{ - \frac{1}{2} + s }
   \left( E_j( R+1, \infty ),\,\mu \right)}\\
\le 
 & \, \widehat{c}_4( s, a_{\rm min}, R, K_{\pm} ) \cdot 
   \Big\{ \| u \|^2_{L^2_{ - \frac{1}{2} - s' }
   ( E_j(R-1,\infty ),\,\mu ) } + \| f \|^2_{ L^2_{ \frac{1}{2} + s }
   ( E_j(R-1,\infty),\,\mu )} \Big\},
\end{align*}
where $\widehat{c}_4(s, a_{\rm min}, R, K_{\pm})$ is a constant depending only on $s,~a_{\rm min},~R$, and $K_{\pm}$. 
\end{prop}
\begin{proof}
Let $\varphi_R(r)$ be the function of $r\ge 0$ defined by
\begin{align*}
 \varphi_R(r) := 
 \begin{cases}
  \ \ \ 0     \qquad & {\rm if} \quad r \le R,\\
  \ r - R \qquad & {\rm if} \quad R \le r \le R + 1,\\
  \ \ \ 1     \qquad & {\rm if} \quad R + 1 \le r.
 \end{cases}
\end{align*}
In Proposition $5.1$, we shall substitute 
\begin{align*}
  \eta (x,z) = p_{\pm}(x,z);~~\psi = \varphi_R (r)\,r^{2s} \exp( -2{\rm Re}\, p_{\pm});~~\Omega = E_j(R,t);~~X = \nabla r,
\end{align*}
where $t > R+1$ is a constant; recall that $p_{\pm}$ is defined by $(11)$. 
For simplicity, we denote $Y := e^{-\eta}\nabla v = \nabla u + u\nabla p_{\pm}$. 
Then, $e^{ -2{\rm Re}\,\eta } = e^{- \eta - \overline{\eta}} = |e^{-\eta}|^2$\,; $\nabla \psi = |e^{-\eta}|^2 \,r^{2s} \bigl\{ \varphi_R' + ( 2sr^{-1}-2 {\rm Re}\, \partial_r \eta ) \varphi_R \bigr\} \nabla r$. 
Hence, each term, appeared in Proposition $5.1$, is calculated as follows:
\begin{align}
 & \psi \Big\{ {\rm Re} \,\langle \nabla r, \nabla \overline{v} \rangle 
   \langle \nabla v, \nabla r \rangle - \frac{1}{2} 
   |\nabla v|^2 \langle \nabla r, \nabla r \rangle \Big\} \\
=& \varphi_R (r)\,r^{2s} |e^{-\eta}|^2 \Big\{ {\rm Re} \,\langle \nabla r, 
   e^{\overline{\eta}}~\overline{Y} \rangle 
   \langle e^{\eta}Y, \nabla r \rangle - \frac{1}{2} |e^{\eta}Y|^2  \Big\}
   \notag \\
=& \varphi_R (r)\,r^{2s} \Big\{ |\langle Y, \nabla r \rangle|^2 
   - \frac{1}{2}|Y|^2 \Big\} \,; \notag \\ 
 & \Big\{ \langle \nabla r, \nabla w \rangle \psi - \frac{1}{2} \langle 
   \nabla r, \nabla \psi \rangle - \frac{1}{2} \psi \Delta_g r \Big\}
   |\nabla v|^2 \\
=& \Big\{ (\partial_r w) \varphi_R \,r^{2s} |e^{-\eta}|^2 
   - \frac{1}{2} \partial_r \psi - \frac{1}{2} \varphi_R \,r^{2s} 
   |e^{-\eta}|^2 \Delta_g r \Big\} |e^{\eta}Y|^2 \notag \\
=& r^{2s} |Y|^2 \Big\{ \varphi_R \big( \partial_r w - sr^{-1} + 
   {\rm Re}\,\partial_r p_{\pm} - \frac{1}{2} \Delta_g r \big) 
   - \frac{1}{2}\varphi_R' \Big\} \,; \notag \\ 
 & {\rm Re} \, \langle \nabla \psi + 2 \psi \nabla \eta, \nabla v \rangle 
   \langle \nabla r, \nabla \overline{v} \rangle \\
=& \,{\rm Re} \,\bigl\langle (\partial_r \psi) \nabla r 
   + 2 \psi (\partial_r \eta) \nabla r, e^{\eta}Y \bigr\rangle 
   \langle \nabla r, e^{\overline{\eta}}~\overline{Y} \rangle \notag \\
=& \,{\rm Re} \, \{ \partial_r \psi + 2 \psi (\partial_r \eta) \}
   \, |e^{\eta}|^2 \,|\langle \nabla r, Y \rangle|^2 \notag \\
=& r^{2s}\,{\rm Re} \, \left\{ \varphi_R' + ( 2sr^{-1} - 2 {\rm Re}\, 
   \partial_r \eta ) \varphi_R + 2\varphi_R\,\partial_r \eta \right\}
   \,|\langle \nabla r, Y \rangle|^2 \notag \\
=& r^{2s}\,\left\{ \varphi_R'+ 2sr^{-1}\varphi_R \right\} 
   \,|\langle \nabla r, Y \rangle|^2 \,; \notag \\ 
 & \psi \,{\rm Re} \, \langle \nabla_{\nabla v} \nabla r, 
   \nabla \overline{v} \rangle 
=  \varphi_R\, r^{2s} |e^{-\eta}|^2 {\rm Re} 
   \langle \nabla_{e^{\eta}Y}\nabla r, e^{\overline{\eta}}~\overline{Y} \rangle
   \\
=& \varphi_R\,r^{2s}\,{\rm Re} \,(\nabla dr)(Y,\overline{Y}) \,; \notag \\ 
 & - \psi \,{\rm Re} \,( qv + e^{\eta} f )\langle \nabla r, 
   \nabla \overline{v} \rangle 
= - \varphi_R\, r^{2s} |e^{-\eta}|^2 \,{\rm Re} \,(qe^{\eta}u + e^{\eta} f )
   \langle \nabla r, e^{\overline{\eta}} \,\overline{Y} \rangle \\ 
=& - \varphi_R\, r^{2s} \,{\rm Re} \,(qu + f) \langle \nabla r, \notag 
  \overline{Y} \rangle .
\end{align}
Also, we have
\begin{align}
  \varphi_R ' = 
  \begin{cases}
  \ 1 \quad & {\rm on} \quad E_j(R,R+1), \\
  \ 0 \quad & {\rm on} \quad M \backslash E_j(R,R+1). 
  \end{cases}
\end{align}
Substituting $(22)$--$(27)$ into the equation of Proposition $5.1$, we obtain 
\begin{align}
 & \int_{S_j(t)} r^{2s} \Big\{ | \langle Y,\nabla r \rangle |^2
   - \frac{1}{2} |Y|^2 \Big\} \,dA_w \\
 & - \int_{E_j(R,R+1)} r^{2s} \Big\{ | \langle Y, \nabla r\rangle |^2
   - \frac{1}{2} |Y|^2 \Big\} \,d\mu \notag \\
 & + \int_{E_j(R,t)} r^{2s} \varphi_R \, 
   {\rm Re}\, ( qu + f ) \langle \overline{Y}, \nabla r \rangle \, d\mu 
   \notag \\ 
=& \int_{E_j(R,t)} r^{2s} \varphi_R \, 
   \Big\{ \big( \partial_r w - sr^{-1} + {\rm Re}\, \partial_r p_{\pm} 
   - \frac{1}{2} \Delta_g r \big) |Y|^2  \notag \\
 & \hspace{30mm} + 2s r^{-1} | \langle Y, \nabla r \rangle |^2
   + {\rm Re}\, (\nabla dr)( Y, \overline{Y} ) \Big\} \,d\mu; \notag 
\end{align}
note that $\varphi_R(R)=0$, and hence, the boundary integral on $S_j(R)$ vanishes. 

We shall write 
\begin{align*}
  Y = \langle Y, \nabla r \rangle \nabla r \oplus Y^{\perp}, 
  \quad {\rm where}~~\nabla r \perp Y^{\perp}. 
\end{align*}
In the following, we shall bound the integrand of the right hand side of $(28)$ from below. 

First, let us consider the case that $1\le j\le m_0$; then, on $E_j$, 
\begin{align*}
 & w \equiv 0 \,; \quad p_{\pm} = \mp i \sqrt{z}~r + \frac{b_j}{2} \log r(x) \,; \\
 & {\rm Re}\,\partial_r p_{\pm}
   \equiv \pm {\rm Im}\,\sqrt{z} + \frac{b_j}{2r} \ge \frac{b_j}{2r} 
   \quad {\rm for}~~z\in {\it \Pi}_{\pm} \,; \\
 & \Delta_g r = \frac{b_j}{r} + O\left ( r^{-1-\delta} \right) \,; 
  \quad {\rm Re}\,(\nabla dr) (Y, \overline{Y}) 
   \ge
   \big\{ \frac{a_j}{r} + O\left( r^{-1-\delta} \right) \big\} 
   |Y^{\perp}|^2.
\end{align*}
Hence, in this case, the integrand of the right hand side of $(28)$ is bounded from below as follows:
\begin{align}
 & r^{2s} \varphi_R
   \Big\{ \big( \partial_r w - \frac{s}{r} + {\rm Re}\, \partial_r p_{\pm} 
   - \frac{1}{2} \Delta_g r \big)|Y|^2 
   + \frac{2s}{r} | \langle Y, \nabla r \rangle|^2
   + {\rm Re}\,(\nabla dr)(Y,\overline{Y}) \Big\} \\
\ge 
 & r^{2s} \varphi_R 
   \Big\{ \Big( \frac{s}{r} - O\left( r^{-1-\delta} \right) \Big)
   |\langle Y, \nabla r \rangle|^2 
 + \Big( \frac{a_j - s}{r} 
   - O \left( r^{-1-\delta} \right) \Big) |Y^{\perp}|^2 \Big\} \notag \\
\ge 
 & \frac{1}{2} r^{2s-1} \varphi_R \cdot \min\{s,a_j - s\} \cdot |Y|^2 
   \qquad {\rm for}~x \in E_j(R_j', \infty) \ \ ( 1 \le j \le m_0 ), 
   \notag 
\end{align}
where $R_j'>0$ is a constant depending only on the geometry of $E_j$ $(1 \le j \le m_0)$. 

Next, consider the case that $m_0 + 1 \le j \le m$; then, on $E_j$, 
\begin{align*}
 & \partial_rw \equiv \frac{\beta_j}{2} \,; 
   \ p_{\pm} = \mp i \sqrt{z - \frac{(\beta_j)^2}{4} }~r, \ 
  {\rm Re}\,\partial_r p_{\pm} \equiv \pm {\rm Im}\,\sqrt{z - \frac{(\beta_j)^2}{4}} \ge 0 
   \ \ {\rm for}~~z \in {\it \Pi}_{\pm} \,; \\
 & \Delta_g r = \beta_j + O\left( r^{-1-\delta}\right ) \,; 
   \quad {\rm Re}\,(\nabla dr) ( Y, \overline{Y} )
   \ge \left\{ \alpha_i + O\left(r^{-1-\delta} \right) \right\} |Y^{\perp}|^2.
\end{align*}
Hence, in this case, the integrand of the right hand side of $(28)$ is bounded from below as follows:
\begin{align}
 & r^{2s} \varphi_R
   \Big\{ \Big( \partial_r w - \frac{s}{r} 
   + {\rm Re}\, \partial_r p_{\pm} - \frac{1}{2} \Delta_g r \Big) |Y|^2 
   + 2 \frac{s}{r} |\langle Y, \nabla r \rangle|^2
   + {\rm Re}\,(\nabla dr) (Y, \overline{Y}) \Big\} \\
\ge
 & r^{2s} \varphi_R \Big\{ \Big( \frac{\beta_j}{2} - \frac{s}{r} 
   \pm {\rm Im} \sqrt{ z - \frac{(\beta_j)^2}{4} } - \frac{\beta_j}{2} 
   - O\left( r^{-1-\delta} \right) \Big) |Y|^2 \notag \\
 & \hspace{30mm} + 2\frac{s}{r} |\langle Y, \nabla r \rangle|^2
   + \big( \alpha_i - O\left( r^{-1-\delta} \right) \big) |Y^{\perp}|^2 
   \Big\} \notag \\
\ge 
 & \frac{1}{2} sr^{2s-1} \varphi_R |Y|^2 
   \qquad {\rm for}~x \in E_j (R_j, \infty) \ \ ( m_0+1 \le j \le m ), \notag 
\end{align}
where $R_j > 0$ is a constant depending only on the geometry of $E_j$. 
Thus, $(29)$ and $(30)$ imply that 
\begin{align}
& r^{2s} \varphi_R \biggl\{ \left( \partial_r w - \frac{s}{r} 
  + {\rm Re}\, \partial_r p_{\pm} - \frac{1}{2} \Delta_g r \right) |Y|^2 
  + 2 \frac{s}{r} |\langle Y,\nabla r\rangle|^2
  + {\rm Re}\,(\nabla dr) (Y, \overline{Y}) \biggr\} \\
\ge 
& \, \widehat{c}_1(s) \cdot r^{2s-1} \varphi_R |Y|^2 
  \qquad {\rm for~any}~x \in M~{\rm with}~r(x) \ge R_3 > 0, \notag 
\end{align}
where $R_3$ is a constant depending only on the geometry of $M$ and $\widehat{c}_1(s) := \frac{ \min\{ s, a_{\rm min} - s \}}{2}$. 

Now, we shall bound the left hand side of $(28)$ from above. 
We begin with the third term of the left hand side of $(28)$. 
As we have seen in Lemma $5.1$, there exists a constant $\widetilde{c} > 0$ such that $|q| \le \widetilde{c} \,r^{-1-\delta}$ for $r(x)\ge 1$. 
Hence, by Schwarz's inequality, we have
\begin{align*}
 & |q||u||Y| \le \widetilde{c} \,r^{-1-\delta} |u| |Y| 
   \le \frac{\widetilde{c}\,^2}{\widehat{c}_1(s)} r^{-1-2\delta} |u|^2 
   + \frac{\widehat{c}_1(s)}{4} r^{-1}|Y|^2 ;\\
 & |f||Y| \le \frac{r}{\widehat{c}_1(s)} |f|^2 
   + \frac{\widehat{c}_1(s)}{4}r^{-1}|Y|^2 .
\end{align*}
From these inequalities, we have
\begin{align}
 & r^{2s} \varphi_R\,{\rm Re}\,(qu+f) 
   \langle \overline{Y}, \nabla r \rangle \\
\le 
 & r^{2s} \varphi_R \left\{ \frac{\widehat{c}_1(s)}{2} r^{-1} |Y|^2 
   + \widehat{c}_2(s) r^{-1-2\delta} |u|^2 + \widehat{c}_2(s) r|f|^2 \right\}, 
   \nonumber 
\end{align}
where $\widehat{c}_2(s) := \frac{ \max \{ (\widetilde{c})^2 , 1 \}}{ \widehat{c}_1(s)}$. 

Next, we bound the second term of the left hand side of $(28)$. 
Since $Y = e^{-\eta}\nabla v = \nabla u + u\nabla p_{\pm}$, we have 
\begin{align}
  \frac{1}{2} \int_{E_j(R,R+1)} r^{2s} |Y|^2 \,d\mu
  \le 
  \int_{E_j(R,R+1)} r^{2s}
  \bigl\{ |\nabla u|^2 + |\nabla p_{\pm}|^2|u|^2 \bigr\} \,d\mu.
\end{align}
In Lemma $4.1$, we set 
\begin{align*}
  \varphi (r) =
 \begin{cases}
  \ r-R+1  \quad & {\rm if}~~R-1 \le r \le R, \\
  \ \ \ \ \ \ \ 1      \quad & {\rm if}~~R \le r \le R+1, \\
  \ -r+R+2 \quad & {\rm if}~~R+1 \le r \le R+2, \\
  \ \ \ \ \ \ \ 0      \quad & {\rm otherwise},
 \end{cases}
\end{align*}
and substitute $\Omega =E_j$, $\varepsilon = \frac{1}{2}$, and $- Lu = zu + f$; then, we obtain 
\begin{align}
 & \int_{E_j(R, R+1)} r^{2s} |\nabla u|^2 \,d\mu \\
\le 
 & \big\{ |z|^2 + 10 + 2 \max_M |V| \big\} \int_{E_j(R-1, R+2)}   r^{2s} |u|^2 \,d\mu  + \int_{E_j(R-1, R+2)} r^{2s} |f|^2 \,d\mu ,\notag 
\end{align}
where we have used the assumptions, $R\ge 2$ and $0 < s \le 1/2$. 
Combining $(33)$ and $(34)$ yields
\begin{align}
& \frac{1}{2} \int_{E_j(R, R+1)} r^{2s} |Y|^2 \,d\mu \\
\le 
& \widehat{c}_3 \int_{E_j(R-1, R+2)} r^{-1-2s'} |u|^2 \,d\mu
   + \int_{E_j(R-1, R+2)} r^{1+2s} |f|^2 \,d\mu, \notag 
\end{align}
where we have again used the assumption, $R\ge 2$, and set 
\begin{align*}
 \widehat{c}_3 
 := \Big\{ 10 + |z|^2 + 2\max_M |V| + \max_{E_j(R, R+1)} |p_{\pm}|^2 \Big\}
  (R + 2)^{ 1 + 2s + 2s' }.
\end{align*}
Putting together $(28)$, $(31)$, $(32)$, and $(35)$, we obtain 
\begin{align}
 & \frac{1}{2} \widehat{c}_1(s) \int_{E_j(R+1, t)} r^{2s-1} |Y|^2 \,d\mu \\
\le 
 & \widehat{c}_2(s) \int_{E_j(R,t)} r^{ - 1 + 2s - 2\delta } |u|^2 \,d\mu
   + \widehat{c}_2(s) \int_{E_j(R,t)} r^{1+2s} |f|^2 \,d\mu \notag \\
 & + \int_{S_j(t)} r^{2s} |\langle Y,\nabla r\rangle|^2 \,dA_w
   + \frac{1}{2} \int_{E_j(R,R+1)} r^{2s} |Y|^2 \,d\mu \notag \\
\le
 & \left\{ ( \widehat{c}_2(s) + \widehat{c}_3 ) \int_{E_j(R-1,R+2)} 
   + \widehat{c}_2(s) \int_{E_j(R+2,t)} \right\} r^{-1-2s'}|u|^2 \,d\mu 
   \notag \\
 & + \left\{ ( \widehat{c}_2(s) + 1 ) \int_{E_j(R-1,R+2)} 
   + \widehat{c}_2(s) \int_{E_j(R+2,t)} \right\} r^{1+2s}|f|^2 \,d\mu 
   \notag \\
& + \int_{S(t)} r^{2s} | \langle Y, \nabla r \rangle |^2 \,dA_w, \notag 
\end{align}
where we have used the assumption $s + s' \le \delta$. 

Since $u$ satisfies the radiation condition (see $(14)$), $\langle Y, \nabla r\rangle = \left( \partial_r + \partial_rp_{\pm} \right) u\in L^2_{-\frac{1}{2}+s}(M,\mu)$, and hence, there exists a divergent sequence $\{t_i\}_{i=1}^{\infty}$ of positive numbers such that $\lim_{i\to \infty} \int_{S(t_i)} r^{2s} |\langle Y, \nabla r\rangle|^2 \,dA_w = 0$.
Hence, substituting $t=t_i$ in $(36)$ and letting $i\to \infty$, we obtain 
\begin{align*}
 & \frac{\widehat{c}_1(s)}{2} \int_{E_j(R+1,\infty)} r^{2s-1}
   \left| \nabla u + u \nabla p_{\pm} \right|^2 \,d\mu \\
\le 
 & \bigl( \widehat{c}_2(s) + \widehat{c}_3 \bigr) \int_{E_j(R-1,\infty)} 
   r^{-1-2s'} |u|^2 \,d\mu + \left( \widehat{c}_2(s) + 1 \right) \int_{E_j(R-1,\infty)}   r^{1+2s} |f|^2 \,d\mu . 
\end{align*}
Proposition $5.2$ follows from this inequality. 
\end{proof}

\section{Decay estimate for the case ${\rm Re}\,z < \frac{(\beta_j)^2}{4}$}
This section studies the decay estimate of a solution $u$ of $(10)$ on an end $E_j$ satisfying ${\rm MC}(\alpha_j ; \beta_j ; \delta_j)$ in the case that ${\rm Re}\,z < \frac{(\beta_j)^2}{4}$. 

First, let us recall $(9)$, that is, on an end $E_j$ satisfying ${\rm MC}(\alpha_j ; \beta_j ; \delta_j)$, the real valued function $V = \frac{\beta_j}{2} ( \Delta_g r - \frac{\beta_j}{2} )$ converges to the positive constant $\frac{(\beta_j)^2}{4}$ at infinity. 
An a priori decay estimate of solutions of $(10)$ will be obtained by taking the real part of the integral in Proposition $3.1$: 
\begin{prop}
Let $(N, g_{_{N}})$ be an $n$-dimensional Riemannian manifold with compact connected $C^{\infty}$-boundary $\partial N$. 
Assume that the inward normal exponential map $\exp^{\perp}_{\partial N}:N^{+}(\partial N) \to N$ induces a diffeomorphism, where $N^{+}(\partial N)=\{$ {\rm inward normal vectors to} $\partial N \}$. 
Let $\beta > 0$ be a constant and set $r(*) := {\rm dist}_g (\partial N, *)$\,{\rm ;} $\widetilde{V}  := \frac{\beta}{2} ( \Delta_{g_{_{N}}} r - \frac{\beta}{2} )$\,{\rm ;} $N(t,\infty) := \{ x\in N \mid r(x) > t \}$ for $t>0$\,{\rm ;} $\mu := e^{-\beta r} v_{g_{_{N}}}$, where $v_{g_{_{N}}}$ stands for the Riemannian measure of $(N, g_{_{N}})$. 
Let $s$ and $s'$ be constants satisfying $0 < s' \le s < 1 , ~s' + s \le 1$, and $u$ be a solution of the equation $- \widetilde{L} u - z u = f$ on $N \backslash \partial N$, where $\widetilde{L} u := \Delta_{g_{_{N}}} u - \beta \,\partial_r u - \widetilde{V} u$. 
Assume that there exist positive constants $\theta$, $\varepsilon$, and $R_1$ such that 
\begin{align*}
 2 \theta > \beta \varepsilon \,; 
 \ {\rm Re}\,z \le \frac{\beta ^2}{4} - \theta \,; \ 
 \Delta_{g_{_{N}}} r \ge \beta - \varepsilon \ {\rm on} \ N(R_1,\infty).
\end{align*}
Then, for any $\ell \ge 2$ and $R \ge \max\{ R_1,2 \}$, we obtain 
\begin{align*}
  \|u\|_{ L^2_{-\frac{1}{2}-s' } ( N(\ell R,\infty),\mu) }^2
\le 
  \frac{\widehat{c}_0 (\theta, \beta, \varepsilon)}{(\ell -1)R}
  \left\{ \|u \|_{L^2_{-\frac{1}{2}-s'}(N(R,\infty),\mu)}^2
  + \|f\|_{L^2_{\frac{1}{2} + s}(N(R,\infty),\mu)}^2 \right\},
\end{align*}
where $\widehat{c}_0 (\theta, \beta, \varepsilon) := \frac{9}{2\theta - \beta \varepsilon}$. 
\end{prop}
\begin{proof}
First, note the following: 
\begin{align}
  \widetilde{V} - {\rm Re}\,z \ge \frac{\beta}{2} 
  \big( \beta - \varepsilon - \frac{\beta}{2} \big) 
  - \frac{\beta^2}{4} + \theta = \theta - \frac{\beta \varepsilon}{2} 
  >0 \quad {\rm on}~~N(R_1,\infty). 
\end{align}
Let $h=h(r)$ be any real valued function of $r\in (R_1,\infty)$ with compact support in $(R_1,\infty)$. 
We shall set $\varphi=h^2$ and $\Omega =E$ in Proposition $3.1$ and take its real part. 
Then, by $(37)$, we have 
\begin{align}
   \Big( \theta - \frac{\beta \varepsilon}{2} \Big) \int_{N} |hu|^2 \,d\mu
\le 
 & \int_{N} ( \widetilde{V} - {\rm Re}\,z )|hu|^2 \,d\mu \\
=& \int_{N} \Bigl\{ h^2 {\rm Re}\,f \overline{u} - h^2 |\nabla u|^2
   - 2hh' {\rm Re} (\partial_r u) \overline{u} \Bigr\} \,d\mu \notag \\
\le 
 & \int_{N} \Bigl\{ (h')^2 |u|^2 + h^2 {\rm Re}\,f\overline{u} \Bigr\} \,d\mu. 
   \notag 
\end{align}
For any constants $\ell \ge 2$, $R\ge 2$ and $t$ satisfying $R_1<R<\ell R<t-1$, set
\begin{align*}
 h(r)=
 \begin{cases}
  \qquad 0 \quad & {\rm if}~~r \le R, \\
  \  \frac{1}{(\ell -1)R}(r-R) \qquad & {\rm if}~~R \le r\le \ell R, \\
  \qquad 1 \quad & {\rm if}~~\ell R\le r \le t-1, \\
  \  - r + t \quad & {\rm if}~~t-1 \le r \le t, \\
  \qquad 0 \quad & {\rm if}~~t\le r,
 \end{cases}
\end{align*}
and substitute this function $h$ into $(38)$. 
Then, we have 
\begin{align}
 & \Big( \theta - \frac{\beta \varepsilon}{2} \Big) 
   \int_{N(\ell R,t-1)} |u|^2 \,d\mu \\
\le 
 & \frac{1}{(\ell -1)^2R^2} \int_{N(R,\ell R)} |u|^2 \,d\mu
   + \int_{N(t-1, t)} |u|^2 \,d\mu + \int_{N(R, t)} |uf| \,d\mu. \notag 
\end{align}
We shall multiply both sides of $(39)$ by $(1+t)^{-2-2s'}$ and integrate it with respect to $t$ over $[\ell R+1,\infty)$; for convenience, we set $\displaystyle F(r) := \int_{S(r)} |u|^2 e^{-2cr} dA$, where $A$ stands for the Riemannian measures induced on each level hypersurface $S(r)=\{ x \in N \mid {\rm dist}\,(\partial N, x) = r \}$. 
Then, as for the integral of the left hand side of $(39)$, we obtain, by Fubini's theorem, 
\begin{align}
 & \int_{\ell R+1}^{\infty} (1+t)^{-2-2s'} \,dt
   \int_{N(\ell R,t-1)} |u|^2 \,d\mu 
=  \int_{\ell R+1}^{\infty} (1+t)^{-2-2s'} \,dt
   \int_{\ell R}^{t-1} F(r) \,dr \\
=& \int_{\ell R}^{\infty} F(r) \,dr 
   \int_{r+1}^{\infty} (1+t)^{-2-2s'} \,dt 
=  \frac{1}{1+2s'} \int_{\ell R}^{\infty} (2+r)^{-1-2s'} F(r) \,dr
   \notag \\
=& \frac{1}{1+2s'} 
   \int_{N(\ell R, \infty)} (2+r)^{-1-2s'} |u|^2 \,d\mu \notag \\
\ge 
 & \frac{1}{1+2s'}\left(\frac{\ell R+1}{\ell R+2}\right)^{1+2s'}
   \int_{N(\ell R,\infty)} (1+r)^{-1-2s'} |u|^2 \,d\mu \notag \\
\ge 
 & \frac{1}{3} \int_{N(\ell R,\infty)} 
   (1+r)^{-1-2s'} |u|^2 \,d\mu, \notag 
\end{align}
where note that $\ell \ge 2$, $R \ge 2$, and $0 < s' \le \frac{1}{2}$. 

As for the integral of the first term on the right hand side of $(39)$, we have
\begin{align}
 & \frac{1}{(\ell -1)^2R^2} \int_{N(R,\ell R)} |u|^2 \,d\mu
   \int_{\ell R+1}^{\infty} (1+t)^{-2-2s'} \,dt \\
=& \frac{(\ell R+2)^{-1-2s'}}{(1+2s')(\ell -1)^2R^2} 
   \int_{N(R,\ell R)} |u|^2 \,d\mu \notag \\
\le 
 & \frac{1}{(1+2s')(\ell -1)^2R^2} 
   \left( \frac{1+\ell R}{\ell R+2} \right)^{1+2s'} 
   \int_{N(R,\ell R)} (1+r)^{-1-2s'} |u|^2 \,d\mu \notag \\
\le 
 & \frac{1}{(\ell -1)^2R^2} \int_{N(R,\ell R)} 
   (1+r)^{-1-2s'} |u|^2 \,d\mu \notag 
\end{align}
As for the integral of the second term on the right hand side of $(39)$, we have
\begin{align}
 & \int_{\ell R+1}^{\infty} (1+t)^{-2-2s'} \,dt 
   \int_{N(t-1,t)} |u|^2 \,d\mu 
=  \int_{\ell R+1}^{\infty} (1+t)^{-2-2s'} \,dt 
   \int_{t-1}^t F(r) \,dr \\
\le 
 & \int_{\ell R}^{\infty} F(r) \,dr \int_r^{r+1}(1+t)^{-2-2s'} \,dt 
\le 
   \int_{\ell R}^{\infty} (1+r)^{-2-2s'} F(r) \,dr \notag \\
=& \int_{N(\ell R,\infty)} (1+r)^{-2-2s'} |u|^2 \,d\mu 
\le 
   \frac{1}{\ell R+1} \int_{N(\ell R,\infty)} (1+r)^{-1-2s'} 
   |u|^2 \,d\mu \notag \\
\le 
 & \frac{1}{(\ell -1)R}\int_{N(\ell R,\infty)}(1+r)^{-1-2s'}|u|^2d\mu,
   \notag 
\end{align}
where, in the first inequality, we have used Fubini's theorem; in the second inequality, we have used the fact that $\int_{r}^{r+1} (1+t)^{-2-2s'} \,dt \le (1+r)^{-2-2s'} \int_{r}^{r+1}\,dt = (1+r)^{-2-2s'}$. 

As for the integral of the third term on the right hand side of $(39)$, by Fubini's theorem, we have 
\begin{align}
 & \int_{\ell R+1}^{\infty} (1+t)^{-2-2s'} \,dt \int_{N(R,t)} |uf| \,d\mu  
=  \int_{\ell R+1}^{\infty} (1+t)^{-2-2s'} \,dt \int_R^t \widetilde{F}(r) \, dr
   \\
=& \int_R^{\ell R+1} \widetilde{F}(r) \, dr 
   \int_{\ell R+1}^{\infty} (1+t)^{-2-2s'} \,dt 
   + \int_{\ell R+1}^{\infty} \widetilde{F}(r) \, dr 
   \int_r^{\infty} (1+t)^{-2-2s'} \,dt \notag \\
=& \frac{1}{(1+2s')(\ell R+2)^{1+2s'}} \int_R^{\ell R+1} \widetilde{F}(r) \, dr
   + \frac{1}{1+2s'} \int_{\ell R+1}^{\infty} (1+r)^{-1-2s'} \widetilde{F}(r) 
   \, dr \notag \\ 
=& \frac{1}{1+2s'} \left\{ \frac{1}{(\ell R+2)^{1+2s'}} \int_{N(R,\ell R+1)} |uf| \,d\mu 
+ \int_{N(\ell R+1,\infty)} (1+r)^{-1-2s'} |uf| \,d\mu\right\}, \notag 
\end{align}
where we have set $\displaystyle\widetilde{F}(r) := \int_{S(r)} |uf| e^{-2cr} dA$. 
Since $s'\le s$, we have
\begin{align}
  \int_{N(R,\ell R+1)} |uf| \,d\mu 
\le \frac{1}{2}\int_{N(R,\ell R+1)} \left\{ (1+r)^{-1-2s'}|u|^2 
  + (1+r)^{1+2s}|f|^2 \right\} \,d\mu.
\end{align}
Also, we have
\begin{align}
& \int_{N(\ell R+1,\infty)} (1+r)^{-1-2s'} |uf| \,d\mu 
  = \int_{N(\ell R+1,\infty)} (1+r)^{(-1-s'-s)+(s-s')} |uf| \,d\mu \\
\le 
 & \frac{1}{(2+\ell R)^{1+s+s'}} \int_{N(\ell R+1,\infty)} (1+r)^{s-s'} 
   |uf| \,d\mu \notag \\
\le 
 & \frac{1}{2(2+\ell R)^{1+s+s'}} \int_{N(\ell R+1,\infty)} \left\{ 
   (1+r)^{-1-2s'}|u|^2 + (1+r)^{1+2s}|f|^2 \right\} \,d\mu. \notag 
\end{align}
Therefore, since $s'\le s$, by combining $(43)$, $(44)$, and $(45)$, we obtain 
\begin{align}
 & \int_{\ell R+1}^{\infty} (1+t)^{-2-2s'} \,dt \int_{N(R, t)} |uf| \,d\mu \\
\le 
 & \frac{1}{2(1+2s')(\ell R+2)^{1+2s'}} \int_{N(R,\infty)} \left\{ 
   (1+r)^{-1-2s'} |u|^2 + (1+r)^{1+2s} |f|^2 \right\} \,d\mu \notag \\
\le 
 & \frac{1}{2(\ell -1)R} \int_{N(R,\infty)} \left\{ (1+r)^{-1-2s'} |u|^2 
   + (1+r)^{1+2s} |f|^2 \right\} \,d\mu. \notag 
\end{align}
Now, putting together $(40),(41),(42)$ and $(46)$, we obtain 
\begin{align*}
 & \frac{1}{3} \Big( \theta - \frac{\beta \varepsilon}{2} \Big) \int_{N(\ell R,\infty)} (1+r)^{-1-2s'} 
   |u|^2 \,d\mu \\
\le 
 & \frac{3}{2(\ell -1)R} \int_{N(R, \infty)} (1+r)^{-1-2s'} |u|^2 \,d\mu 
  + \frac{1}{2(\ell -1)R} \int_{N(R, \infty)} (1+r)^{1+2s} |f|^2 \,d\mu.
\end{align*}
Proposition 6.1 follows from this inequality. 
\end{proof}

\section{Decay estimate for the case that ${\rm Re} \,z > \frac{(\beta_j)^2}{4}$}
In this section, we will show an a-priori-decay-estimate of an outgoing or incoming solution of the equation $(10)$ on the end $E_j$ in the case that ${\rm Re} \,z > \frac{(\beta_j)^2}{4}$. 
This will be accomplished by combining Proposition $5.2$ and Lemma $7.1$ below; lemma $7.1$ will be proved by taking the imaginary part of the integral in Proposition $3.1$. 
\begin{prop}
Let $2 \le R$, $z \in K_{\pm}$, and $f\in L^2_{\frac{1}{2}+s}(M,\mu)$, and $u$ be a solution of $(10)$ satisfying the radiation condition. 
Assume that 
\begin{align*}
  S < a_{\min}; \quad s + s' \le \delta; \quad 
{\rm Re} \,z > \frac{(\beta_j)^2}{4}.
\end{align*}
Then, for any $R$ and $R_1$ satisfying $\max\{1, r_0\} \le R_1<R_1+1<R$, we obtain 
\begin{align*}
& {\rm Re} \,\sqrt{ z - \frac{(\beta_j)^2}{4} } 
  ~\|u\|^2_{L^2_{-\frac{1}{2}-s'}
  (E_j(R,\infty),\mu)} \\ 
\le 
& \, \widehat{c}_3 \,(1+R)^{-2s'} \left\{ \|u\|^2_{L^2_{-\frac{1}{2}-s'}
  (E_j(R_1-1,\infty),\mu)} + \|f\|^2_{L^2_{\frac{1}{2}+s}
  (E_j(R_1-1,\infty),\mu)} \right\},
\end{align*}
where $\widehat{c}_3 = \widehat{c}_3 (s, s', a_{\min}, R_1,K_{\pm})$ is a constant depending only on $s$, $s'$, $a_{\min}$, $R_1$, and $K_{\pm}$. 
\end{prop}
As is mentioned above, Proposition $7.1$ immediately follows from the following Lemma $7.1$ and Proposition $5.2$, and hence, it suffices to prove the following:
\begin{lem}
Let $f \in L^2_{\frac{1}{2}+s}(M,\mu)$, $z \in {\it \Pi}_{\pm}$, and ${\rm Re}\, z > \frac{(\beta_j)^2}{4}$. 
Assume that $u$ is an outgoing or incoming solution of $(10)$. 
Then, for any $R$ and $R_1$ satisfying $\max\{1, r_0\} \le R_1 < R_1+1 < R$, we obtain 
\begin{align*}
 & {\rm Re} \sqrt{ z - \frac{(\beta_j)^2}{4} }~\|u\|^2_{L^2_{-\frac{1}{2}-s'}(E_j(R,\infty),\mu)}
   \\
\le 
 & \frac{(1+R)^{-2s'}}{2s'} \left\{ \| (\partial_r + \partial_r p_{\pm})u \|^2
   _{L^2_{-\frac{1}{2}+s}(E_j(R,\infty),\mu)} 
   + \|f\|^2_{L^2_{\frac{1}{2}+s}(E_j(R_1-1,\infty),\mu)}\right.\\
 & \hspace{75mm}\left. + \widehat{c}_4 \|u\|^2_{L^2_{-\frac{1}{2}-s'}(E_j(R_1-1,\infty),\mu)} \right\},
\end{align*}
where $\widehat{c}_4 := 64 \,(R_1)^2 \left\{ |z|^2 + 1 + \max_M |V|\right\}$. 
\end{lem}
\begin{proof}
Let $t$ be a constant satisfying $t > R_1+1$. 
In proposition $3.1$, we shall set $\Omega = E_j(R_1,t)$ and 
\begin{align*}
 \varphi(r) = 
 \begin{cases}
  \ \ \ 0 \quad & {\rm if}~~R_1 \le r,\\
    r - R_1 \quad & {\rm if}~~R_1 \le r \le R_1 + 1,\\
  \ \ \ 1 \quad & {\rm if}~~R_1 + 1 \le r ,
 \end{cases}
\end{align*}
and take the imaginary part. 
Then, as in the proof of Proposition $3.2$, we obtain 
\begin{align}
 & \int_{S_j(t)} {\rm Im} \,(\partial_r p_{\pm}) |u|^2 \,dA_w 
   - {\rm Im} \,z \int_{E_j(R_1,t)} \varphi |u|^2 \,d\mu \\
=& \int_{S_j(t)} {\rm Im} \,\overline{u} \,(\partial_r + \partial_r p_{\pm}) u 
   \,dA_w + \int_{E_j(R_1,t)} \varphi \,{\rm Im} \,f \overline{u} \,d\mu 
   - \int_{E_j(R_1,R_1+1)} {\rm Im} \,(\partial_r u) \overline{u} \,d\mu. 
   \nonumber 
\end{align}
Here, note that 
\begin{align}
& {\rm Im} \,z \ge 0, \quad {\rm Im} \,\partial_r p_{+} 
  = - {\rm Re} \,\sqrt{ z - \frac{(\beta_j)^2}{4} } < 0 
  \qquad {\rm if}~~z \in {\it \Pi}_{+}; \\
& {\rm Im} \,z \le 0, \quad {\rm Im} \,\partial_r p_{-} 
  = {\rm Re} \,\sqrt{ z - \frac{(\beta_j)^2}{4} } > 0 
  \qquad {\rm if}~~z\in {\it \Pi}_{-};
\end{align}
$(48)$ and $(49)$ follow from the assumption ${\rm Re}\,z > \frac{(\beta_j)^2}{4}$. 
Thus, signs of real numbers, ${\rm Im}\,z$ and ${\rm Im} \,\partial_r p_{\pm}$, are different in both cases. 
Hence, by $(47)$, $(48)$, and $(49)$, we obtain 
\begin{align}
 & {\rm Re} \,\sqrt{z-\frac{(\beta_j)^2}{4}} \int_{S_j(t)} |u|^2 \,dA_w \\
\le
 & {\rm Re} \,\sqrt{z-\frac{(\beta_j)^2}{4}} \int_{S_j(t)} |u|^2 \,dA_w 
   + |{\rm Im}\, z| \int_{E_j(R_1,t)} \varphi |u|^2 \,d\mu \notag \\
\le 
 & \int_{S_j(t)} \bigl| \overline{u} \,( \partial_r + \partial_r p_{\pm} ) u 
   \bigr| \,dA_w + \int_{E_j(R_1,t)} |f \overline{u}| \,d\mu 
   + \int_{E_j(R_1,R_1+1)} \bigl| (\partial_r u) \overline{u} \bigr| \,d\mu. 
   \notag 
\end{align}
We shall multiply both sides of the inequality $(50)$ by $(1+t)^{-1-2s'}$ and integrate it over $[R,\infty)$ with respect to $t$. 
Then, as for the integral of the first term on the right hand side of $(50)$, we have, by Schwarz's inequality, 
\begin{align}
& \int_R^{\infty} (1+t)^{-1-2s'} \,dt \int_{S_j(t)} 
  \bigl| u \bigr| \bigl| ( \partial_r + \partial_r p_{\pm} ) u \bigr| \,dA_w \\
\le 
& (1+R)^{-s'-s} \int_R^{\infty} dt \int_{S_j(t)} (1+t)^{-1+s-s'} 
   \bigl| u \bigr| \bigl| ( \partial_r + \partial_r p_{\pm} ) u \bigr| \,dA_w 
   \notag \\ 
\le 
& \frac{(1+R)^{-s'-s}}{2} \biggl\{ \int_{E_j(R,\infty)} (1+r)^{-1+2s} 
  \bigl| ( \partial_r + \partial_r p_{\pm} ) u \bigr|^2 \, d\mu \notag \\
& \hspace{60mm} + \int_{E_j(R,\infty)} (1+r)^{-1-2s'} |u|^2 \, d\mu \biggr\} .
   \notag 
\end{align}
We shall set $\displaystyle \widehat{F}(r) := \int_{S_j(t)} |uf| \,d\mu$ and use Fubini's theorem for the integral of the second term on the right hand side of $(50)$. 
Then, we have 
\begin{align}
 & \int_R^{\infty} (1+t)^{-1-2s'} \,dt \int_{E_j(R_1,t)} |uf| \,d\mu 
 = \int_R^{\infty} (1+t)^{-1-2s'} \,dt \int_{R_1}^{t} \widehat{F}(r) \,dr \\
=& \int_{R_1}^R \widehat{F}(r) \,dr \int_R^{\infty} (1+t)^{-1-2s'} \,dt 
   + \int_{R}^{\infty} \widehat{F}(r) \,dr \int_r^{\infty} (1+t)^{-1-2s'} \,dt   \notag \\
=& \frac{(1+R)^{-2s'}}{2s'} \int_{E_j(R_1,R)} |uf| \,d\mu 
   + \frac{1}{2s'} \int_{E_j(R,\infty)} (1+r)^{-2s'} |uf| \,d\mu \notag \\
\le 
 & \frac{(1+R)^{-2s'}}{2s'} \int_{E_j(R_1,\infty)} |uf| \,d\mu \notag \\
\le 
 & \frac{(1+R)^{-2s'}}{4s'} \left\{ \int_{E_j(R_1,\infty)} (1+r)^{1+2s} |f|^2 
   \,d\mu + \int_{E_j(R_1,\infty)} (1+r)^{-1-2s'} |u|^2 \,d\mu \right\} ,
   \notag 
\end{align}
where, in the last line, we have used Schwarz's inequality and the assumption $s'\le s$. 
As for the integration of the third term on the right hand side of $(50)$, we have 
\begin{align}
 & \int_R^{\infty} (1+t)^{-1-2s'} \,dt \int_{E_j(R_1,R_1+1)} 
   \bigl| (\partial_r u) u \bigr| \,d\mu \\
\le 
 & \frac{(1+R)^{-2s'}}{4s'} \left\{ \int_{E_j(R_1,R_1+1)} |\partial_r u|^2 
   \,d\mu + \int_{E_j(R_1,R_1+1)} |u|^2 \,d\mu \right\} . \notag 
\end{align}
Now, we shall apply Lemma $4.1$ to the first term on the right hand side of $(53)$; set $\Omega = E_j$, $\varepsilon = \frac{1}{2}$, and 
\begin{align*}
 \varphi (r) =
 \begin{cases}
 \ \ \ \ \ 0 \quad & {\rm if}\quad 0 \le r \le R_1-1, \\
   r-R_1+1 \quad & {\rm if}\quad R_1-1 \le r \le R_1, \\
 \ \ \ \ \ 1 \quad & {\rm if}\quad R_1 \le r \le R_1+1, \\
   -r+R_1+2 \quad & {\rm if}\quad R_1+1 \le r \le R_1+2, \\
 \ \ \ \ \ 0 \quad & {\rm if}\quad R_1+2 \le r.
 \end{cases}
\end{align*}
Then, in view of $-L_{\rm loc}u = z u + f$, we obtain 
\begin{align}
 & \int_{E_j(R_1, R_1+1)} |\nabla u|^2 \,d\mu \\
\le 
 & \int_{E_j(R_1-1, R_1+2)} (1+r)^{1+2s} |f|^2 \,d\mu \notag \\
 & \hspace{10mm} + (3+R_1)^{1+2s'} \bigl\{ |z|^2+\widehat{c}_5 \bigr\} 
   \int_{E_j(R_1-1, R_1+2)} (1+r)^{-1-2s'} |u|^2 \,d\mu. \notag 
\end{align}
Here, $\widehat{c}_5 := \frac{1}{2} + 2 \max_M |V| + 4(1+s) \le 7 + 2 \max_M |V|$ by the fact $0<s<1/2$. 
Hence, $(53)$ and $(54)$ imply that
\begin{align}
& \int_R^{\infty} (1+t)^{-1-2s'} \,dt \int_{E_j(R_1,R_1+1)} 
  \bigl| (\partial_r u) u \bigr| \,d\mu \\ 
\le 
& \frac{(1+R)^{-2s'}}{4s'} \bigg\{ \int_{E_j(R_1-1,R_1+2)} (1+r)^{1+2s} |f|^2 
  \,d\mu \notag \\
& \hspace{40mm} + \widehat{c}_6 \int_{E_j(R_1-1,R_1+2)} 
   (1+r)^{-1-2s'} |u|^2 \,d\mu \bigg\} , \notag
\end{align}
where $\widehat{c}_6 := (3+R_1)^{1+2s'} \bigl\{ |z|^2 + \widehat{c}_5 + 1 \bigr\} \le (3+R_1)^2 \left\{ |z|^2 + 8 + 2 \max_M |V| \right\}$. 

Thus, putting together $(50)$, $(51)$, $(52)$, and $(55)$, we obtain 
\begin{align}
 & {\rm Re} \,\sqrt{z - \frac{(\beta_j)^2}{4}} \int_{E_j(R,\infty)} (1+r)^{-1-2s'} |u|^2 \,d\mu    \\
\le
 & \frac{(1+R)^{-s'-s}}{2} \int_{E_j(R,\infty)} (1+r)^{-1+2s} 
   \bigl| ( \partial_r + \partial_r p_{\pm} ) u \bigr|^2 \, d\mu \notag \\
 & \hspace{20mm} + \frac{(1+R)^{-2s'}}{2s'} \int_{E_j(R_1-1,\infty)} 
   (1+r)^{1+2s} |f|^2 \,d\mu    \notag \\
 & \hspace{40mm} + \widehat{c}_7(R) \int_{E_j(R_1-1,\infty)} (1+r)^{-1-2s'} 
   |u|^2 \,d\mu, \notag 
\end{align}
where $\widehat{c}_7(R) := \frac{(1+R)^{-s'-s}}{2} + (1 + \widehat{c}_6) \frac{(1+R)^{-2s'}}{4s'}$. 
Since $s' \le s$ and $0<s'\le \frac{1}{2}$, we see that 
\begin{align}
\frac{(1+R)^{-s'-s}}{2} & \le \frac{(1+R)^{-2s'}}{4s'} \,;\\
2s'(1+R)^{2s'} \cdot \widehat{c}_7(R) & \le 1 + \frac{(3+R_1)^{1+2s'}}{2} \left\{ |z|^2 + 8 + \max_M |V| \right\} .
\end{align} 
Lemma $7.1$ follows from $(56)$, $(57)$, and $(58)$. 
\end{proof}

\section{The proof of main theorems}
In this section, we shall prove main theorems after preparing some lemmas. 

In the following, we shall fix $s', s \in \mathbb{R}$ satisfying 
\begin{align}
  0 < s' < s < \min \left\{ a_{\min}, \frac{1}{2} \right\}
  ; \quad s' + s \le \delta,
\end{align}
and use these numbers in the definition of the radiation condition $(13)$ and $(14)$. 
We denote by $R_{-L}(z)$ the resolvent $(-L-z)^{-1}$ of $-L$ for $z\in {\rm resolv}(-L)$, where ${\rm resolv}(-L)$ stands for the resolvent set of $-L$. 
Moreover, we denote by $\mathbb{B}_{\mu} \left( \frac{1}{2}+s,-\frac{1}{2}-s' \right)$ the space of all bounded linear operators $T$ from $L^2_{\frac{1}{2}+s}(M,\mu)$ into $L^2_{-\frac{1}{2}-s'}(M,\mu)$ with ${\rm Dom}(T) = L^2_{\frac{1}{2}+s}(M,\mu)$. 
$\mathbb{B}_{\mu} \left( \frac{1}{2}+s,-\frac{1}{2}-s' \right)$ is a Banach space by the operator norm $\| * \|_{\mathbb{B}_{\mu} \left( \frac{1}{2}+s,-\frac{1}{2}-s' \right)}$.
\vspace{2mm}

We shall begin with the following:
\begin{prop}
Let $z\in \mathbb{C} \backslash \mathbb{R}$. 
Then, {\rm (i)} and {\rm (ii)} below hold: 
\begin{enumerate}[{\rm (i)}]
  \item Let $f \in L^2(M,\mu)$. 
Assume that $u$ is an outgoing or incoming solution of $(10)$, i.e., $u \in L^2_{- \frac{1}{2} - s' }(M, \mu)$ and $\partial_r u + (\partial_r p_{\pm}) u \in L^2_{-\frac{1}{2}+s}(M,\mu)$. 
Then, $u \in L^2(M,\mu)$\,{\rm ;} in particular, $|{\rm Im}\,z| \,\|u\|_{L^2(M,\mu)} \le \|f\|_{L^2(M,\mu)}$, and $u$ coincides with the $L^2(M,\mu)$-solution $R_{-L}(z)f$. 
  \item Let $f \in L^2_{\gamma}(M,\mu)$ for $\gamma \in (0,1]$. 
Assume that $u$ is an outgoing $($or incoming$)$ solution of $(10)$. 
Then, $u \in L^2_{\gamma}(M,\mu)$ and 
\begin{align*}
  |{\rm Im}\,z| \,\|u\|_{L^2_{\gamma}(M,\mu)} 
\le 
  (1 + \gamma) \, \|f\|_{L^2_{\gamma}(M,\mu)} 
  + 10 \gamma \, \Big\{ 1 + |z| + \max_M \sqrt{|V|} \Big\} \|u\|_{L^2(M,\mu)}.
\end{align*}
\end{enumerate}
\end{prop}
\begin{proof}
(i) We may assume that $u$ is non-trivial. 
Setting $\varphi=1$ in Proposition $3.2$, we obtain
\begin{align*}
  |{\rm Im}\,z| \int_{U(R)} |u|^2 \,d\mu 
\le 
  \int_{S(R)} 
  \big| (\partial_r + \partial_r p_{\pm})u \bigr|\,\bigl| u \big| \,dA_w 
  + \int_{U(R)} |u||f| \,d\mu.
\end{align*}
Here, in view of $s' \le s$, we obtain $\big| (\partial_r + \partial_r p_{\pm})u \bigr|\,\bigl| u \big| \le \frac{1}{2} \left\{ r^{2s} \bigl| (\partial_r + \partial_r p_{\pm})u \bigr|^2 + r^{-2s'} |u|^2 \right\}$; also, $\big\| |u||f| \big\|_{L^2(U(R),\mu)}^2 \le \| u \|_{L^2(U(R),\mu)} \|f\|_{L^2(U(R),\mu)}$. 
Hence, we obtain 
\begin{align}
 & |{\rm Im}\,z| \, \| u \|_{L^2(U(R),\mu)} \\
\le 
 & \frac{1}{2} \| u \|_{L^2(U(R),\mu)}^{-1} 
   \int_{S(R)} \Big\{ r^{2s} \big| (\partial_r 
   + \partial_r p_{\pm})u \big|^2 + r^{-2s'} |u|^2 \Big\}\, dA_w 
  + \|f\|_{L^2(U(R),\mu)}. \notag 
\end{align}
Since $u$ satisfies $(14)$, we see that $\displaystyle\liminf_{R\to \infty} \int_{S(R)} \Big\{ r^{-2s'} \bigl|u\bigr|^2 + r^{2s} \bigl| (\partial_r + \partial_r p_{\pm})u \bigr|^2 \Big\} \,dA_w=0$. 
Therefore, substituting an appropriate divergent sequence $\{ R_i \}_{i=1}^{\infty}$ of positive numbers into $(60)$, we get our desired result. 
\vspace{2mm}

\noindent (ii) The proof of (ii) will be completed through two steps. 
First, we shall show that $u\in L^2_{\frac{\gamma}{2}}(M,\mu)$. 
Note that $-Lu = f + z u \in L^2(M, \mu)$ by (i), which implies that $|\nabla u| \in L^2(M,\mu)$ by Corollary $4.1$. 
Hence, 
$\bigl| (\partial_r + \partial_r p_{\pm})u \bigr| \in L^2(M,\mu)$ by $(11)$. 
Therefore, we see that 
\begin{align}
  \liminf_{R\to \infty} \int_{S(R)} (1+r)^{\gamma} \big| (\partial_r 
  + \partial_r p_{\pm})u \big|\, \big| u \big| \, dA_w = 0 
  \quad {\rm for}~~\gamma \in (0,1].
\end{align}
Now, for $\gamma \in (0,1]$, set $\varphi(r)=(1+r)^{\gamma}$ in Proposition $3.2$. 
Then, 
\begin{align}
 & |{\rm Im}\,z| \, \| u \|_{L_{\frac{\gamma}{2}}^2\left( U(R), \mu \right)} ^2 \\
\le 
 & \int_{S(R)} (1+r)^{\gamma} \big| (\partial_r + \partial_r p_{\pm})u \big|\,   \big| u \big| \, dA_w \notag \\
 & + \int_{U(R)} \Big\{ \gamma (1+r)^{\gamma-1} |\partial_r u|\,|u| 
   + (1+r)^{\gamma} |f||u| \Big\} \,d\mu . \notag
\end{align}
Here, $\int_{U(R)} (1+r)^{\gamma-1} |\partial_r u|\,|u| \,d\mu \le \| \partial_r u \|_{L^2_{-1+\frac{\gamma}{2}}\left( U(R), \mu \right)} \cdot \| u \|_{L^2_{\frac{\gamma}{2}}\left( U(R), \mu \right)}$ and $\int_{U(R)} (1+r)^{\gamma} |f| |u| \,d\mu \le \| f \|_{L^2_{\frac{\gamma}{2}}\left( U(R), \mu \right)} \cdot \| u \|_{L^2_{\frac{\gamma}{2}}\left( U(R), \mu \right)}$. 
Therefore, bearing $(61)$ in mind, substituting an appropriate divergent sequence $\{R_i\}_{i=1}^{\infty}$ into $(62)$, and letting $i\to \infty$, we obtain
\begin{align*}
   |{\rm Im}\,z|\, \|u\|_{L^2_{\frac{\gamma}{2}}(M,\mu)} 
 \le 
   \gamma \, \|\partial_r u\|_{L^2_{-1+\frac{\gamma}{2}}(M,\mu)} 
   + \|f\|_{L^2_{\frac{\gamma}{2}}(M,\mu)} < \infty, 
\end{align*}
where, the last inequality follows from $-1+\frac{\gamma}{2} \le 0$ and $|\nabla u| \in L^2(M,\mu)$. 
Thus, $u\in L^2_{\frac{\gamma}{2}}(M, \mu)$. 

Next, we shall prove that $u\in L^2_{\gamma}(M,\mu)$. 
Note that $-Lu = f + z u \in L^2_{\frac{\gamma}{2}}(M,\mu)$ by the fact, $u\in L^2_{\frac{\gamma}{2}}(M,\mu)$, as is shown above. 
Hence, Corollary $4.1$ implies that $\bigl| (\partial_r + \partial_r p_{\pm})u \bigr| \in L^2_{\frac{\gamma}{2}}(M,\mu)$, which, together with $0 < \gamma \le 1$ and $u\in L^2_{\frac{\gamma}{2}}(M,\mu)$, implies that
\begin{align}
  \liminf_{R\to \infty} \int_{S(R)} (1+r)^{2\gamma} 
  \big| (\partial_r + \partial_r p_{\pm})u \big|\,\big| u \big| \, dA_w = 0 
  \quad {\rm for}~0 < \gamma \le 1.
\end{align}
Now, we shall set $\varphi(r) = (1+r)^{2\gamma}$ in Proposition $3.2$ and repeat the same arguments as above: then, we obtain 
\begin{align}
  |{\rm Im}\,z| \, \| u \|_{L^2_{\gamma}(U(R),\mu)} 
\le 
& \| u \|_{L^2_{\gamma}(U(R),\mu)}^{-1} 
  \int_{S(R)} (1+r)^{2\gamma} \bigl| (\partial_r + \partial_r p_{\pm})u \bigr|
  \, \big| u \big| \, dA_w \\
& + 2 \gamma \, \| \partial_r u \|_{L^2_{-1+\gamma}(U(R),\mu)} + \| f \|_{L^2_{\gamma}(U(R),\mu)}. 
   \notag 
\end{align}
Bearing $(63)$ in mind, substituting an appropriate divergent sequence $\{R_i\}_{i=1}^{\infty}$ into $(64)$, and letting $i\to \infty$, we obtain 
\begin{align}
  |{\rm Im}\,z|\, \|u\|_{L^2_{\gamma}(M,\mu)} 
 \le 
  2\gamma \, \|\partial_r u\|_{L^2_{-1+\gamma}(M,\mu)} 
  + \|f\|_{L^2_{\gamma}(M,\mu)}.
\end{align}
Again, since $-1 + \gamma \le 0$, Corollary $4.1$ implies that 
\begin{align}
  \|\nabla u\|_{L^2_{-1+\gamma}(M,\mu)} 
\le 
& \frac{1}{2} \| f + zu \|_{L^2_{-1+\gamma}(M,\mu)} 
  + \widehat{c}_{8} \,\|u\|_{L^2_{-1+\gamma}(M,\mu)} \\
\le 
& \frac{1}{2} \| f \|_{L^2(M,\mu)} 
  + \Big( \widehat{c}_{8} + \frac{|z|}{2} \Big) \|u\|_{L^2(M,\mu)}, \notag 
\end{align}
where $\widehat{c}_{8}:= \sqrt{\widehat{c}\left( \frac{1}{4} \right)} \le \frac{9}{2} + \frac{3}{2} \max_M \sqrt{|V|}$. 
Combining $(65)$ and $(66)$, we obtain
\begin{align*}
& |{\rm Im}\,z|\, \|u\|_{L^2_{\gamma}(M,\mu)} \\
\le 
& \gamma \, \| f \|_{L^2(M,\mu)} 
  + 10 \gamma \, \Big\{ 1 + |z| + \max_M \sqrt{|V|} \, \Big\} \| u \|_{L^2(M,\mu)}   + \| f \|_{L^2_{\gamma}(M,\mu)}. 
\end{align*}
Assertion (ii) follows form this inequality. 
\end{proof}

Proposition $8.1$ implies the following two corollaries:
\begin{cor}
Assume that $z \in {\it \Pi}_{\pm} \backslash (0,\infty)$. 
Then, for any $f \in L^2_{\frac{1}{2}+s}(M,\mu)$, the equation $(10)$ has a unique outgoing or incoming solution $u$\,{\rm ;} moreover, $u$ coincides with the $L^2(M,\mu)$-solution $R_{-L}(z)f$ and belongs to $L^2_{\frac{1}{2}+s}(M,\mu)$. 
\end{cor}
\begin{proof}
Let $f\in L^2_{\frac{1}{2}+s}(M,\mu)$ and consider $L^2(M,\mu)$-solution $R_{-L}(z)f$. 
By Proposition $8.1$, it suffices to prove that $R_{-L}(z)f$ satisfies the radiation condition. 
Then, $R_{-L}(z)f$ satisfies $-L\bigl(R_{-L}(z)f\bigr) = z\,\bigl(R_{-L}(z)f\bigr) + f \in L^2(M,\mu)$, and hence, Corollary $4.1$ implies that $\left| \nabla \big( R_{-L}(z)f \big) \right| \in L^2(M,\mu)$. 
In view of the facts, $0 < s \le \frac{1}{2}$ and $\sup_{x\in M} |\partial_r p_{\pm}(x, z)| < \infty$ for each fixed $z$, we see that $\partial_r \bigl(R_{-L}(z)f\bigr) + ( \partial_r p_{\pm}(z,*)) \big( R_{-L}(z)f \big) \in L^2(M,\mu) \subset L^2_{-\frac{1}{2}+s}(M,\mu)$. 
Thus, $R_{-L}(z) f$ satisfies the condition $(14)$, and hence, $R_{-L}(z) f \in L^2_{\frac{1}{2}+s}(M,\mu)$ by Proposition $8.1$. 
\end{proof}
\begin{cor}[{\bf decay estimate for} $|{\rm Im}\,z|>0$]
Let $z \in {\it \Pi}_{\pm} \backslash (0,\infty)$, $f \in L^2_{\frac{1}{2}+s}(M, \mu)$, and $u$ be the corresponding outgoing or incoming solution of $(10)$. 
Then, for any $R>1$, 
\begin{align*}
  |{\rm Im}\,z|^2 \int_{E(R,\infty)} (1+r)^{-1-2s'} |u|^2 \,d\mu 
\le 
  (\widehat{c}_{11})^2 (1+R)^{-2-2s-2s'} \int_M (1+r)^{1+2s} |f|^2 \,d\mu,
\end{align*}
where $\widehat{c}_{11} := 10 \left\{ 1 + \frac{1}{|{\rm Im}\,z|} \left( 1 + |z| + \max_M |V| \right) \right\}$. 
\end{cor}
\begin{proof}
Corollary $8.1$ and Proposition $8.1$ imply that 
\begin{align}
   |{\rm Im}\,z|\, \|u\|_{L^2_{\frac{1}{2}+s}(M,\mu)} 
\le 
 & 10 \bigg\{ \|f\|_{L^2_{\frac{1}{2}+s}(M,\mu)} 
  + \frac{1 + |z| + \displaystyle \max_M \sqrt{|V|}}{|{\rm Im}\,z|} \|f\|_{L^2(M,\mu)} \bigg\} \\
\le 
 & \, \widehat{c}_{11} \, \| f \|_{L^2_{\frac{1}{2}+s}(M,\mu)}, \notag 
\end{align}
where $\widehat{c}_{11} = 10 \left\{ 1 + \frac{1}{|{\rm Im}\,z|} \left( 1 + |z| + \displaystyle \max_M \sqrt{|V|} \right) \right\}$. 
Therefore, we obtain 
\begin{align*}
 & |{\rm Im}\,z|^2 (1+R)^{2+2s+2s'} \int_{E(R,\infty)} (1+r)^{-1-2s'} |u|^2 
   \,d\mu \\
\le 
 & |{\rm Im}\,z|^2 \int_{E(R,\infty)} (1+r)^{1+2s} |u|^2 \,d\mu 
\le 
  |{\rm Im}\,z|^2 \int_M (1+r)^{1+2s} |u|^2 \,d\mu \\
\le 
 & (\widehat{c}_{11})^2 \int_M (1+r)^{1+2s} |f|^2 \,d\mu 
   \qquad ({\rm by}~(67)).
\end{align*}
Lemma $8.2$ follows from this inequality.
\end{proof}
Next, we shall show the following uniqueness theorem: 
\begin{lem}[{\bf uniqueness}]
Assume that $z = \lambda \in {\it \Pi}_{\pm}\cap (0,\infty)$. 
Then, outgoing or incoming solution $u$ of $(10)$, if it exists, is uniquely determined by $z$ and $f$. 
\end{lem}
\begin{proof}
First, note that $z = \lambda > \frac{(\beta_1)^2}{4}$ by the definition of ${\it \Pi}_{\pm}$. 
Let $u_1$ and $u_2$ be two solutions of the same equation $(10)$. 
Then, $u := u_1 - u_2$ is a solution to the eigenvalue equation 
\begin{align}
  - Lu - \lambda u = 0\,; \qquad \lambda > \frac{(\beta_1)^2}{4}.
\end{align}
Since $\lambda \in \mathbb{R}$, we may assume that $u$ is real-valued by considering the real and imaginary part of $u$. 
Hence, setting $\varphi \equiv 1$ and $f = 0$ in Proposition $3.2$, we obtain, for $R \ge r_0$, 
\begin{align*}
  \int_{S(R)} | {\rm Im} \, \partial_r p_{\pm} |\,|u|^2 \,dA_w 
\le 
  \int_{S(R)} |(\partial_r + \partial_r p_{\pm})u|\,|u| \,dA_w.
\end{align*}
Therefore, in view of $(11)$ and $(12)$, we get, for $R \ge r_0$, 
\begin{align}
  \sqrt{\lambda - \frac{(\beta_{1})^2}{4}} \int_{S_{1}(R)} |u|^2 \,dA_w 
\le 
  \int_{S(R)} |(\partial_r + \partial_r p_{\pm})u|\,|u| \,dA_w.
\end{align}
Since $( \partial_r  + \partial_r p_{\pm} )u \in L^2_{-\frac{1}{2}+s}(M, \mu)$ 
and $u \in L^2_{-\frac{1}{2} - s'}(M, \mu)$ by $(14)$, we see that 
\begin{align*}
 & \int_M (1+r)^{s - s' -1} |(\partial_r + \partial_r p_{\pm}) u|\,|u| 
   \,d\mu \\
\le 
 & \frac{1}{2} \int_M \left\{ (1+r)^{-1 + 2s} |(\partial_r + \partial_r p_{\pm}) u|^2 + (1+r)^{-1 -2s'} |u|^2 \right\} \,d\mu < \infty.
\end{align*}
Hence, multiplying both sides of $(69)$ by $(1+R)^{s-s'-1}$ and integrating it over $[r_0, \infty)$ with respect to $R$, we obtain 
\begin{align*}
& \sqrt{\lambda - \frac{(\beta_{1})^2}{4}} 
  \int_{E_{1}(r_0, \infty)} (1+r)^{s-s'-1} |u|^2 \,d\mu \\ 
\le 
& \int_{E(r_0, \infty)} (1+r)^{s-s'-1} 
  |(\partial_r + \partial_r p_{\pm})u| \,|u| \,d\mu < \infty.
\end{align*}
Thus, $Lu = - \lambda u \in L^2_{\frac{s-s'-1}{2}}(E_{1},\mu)$; hence, $|\nabla u|\in L^2_{\frac{s-s'-1}{2}}(E_{1}, \mu)$ by the same arguments of the proof of Corollary $4.1$. 
Therefore, we obtain
\begin{align}
  \liminf_{R\to \infty} R^{\, s - s'} \int_{S_{1}(R)}
  \bigl\{ |\partial_ru|^2 + |u|^2 \bigr\} \,d\mu = 0,
\end{align}
where $s-s' > 0$ (see $(59)$). 

Now, by reconsidering the arguments developed in \cite{K3} and \cite{K4}, it is not hard to see that, if $(68)$ holds on $E_1$, then $(70)$, together with $(1)$ or $(2)$ with $E=E_1$, implies that $u \equiv 0$ on $E_1$ (in both cases, $(1)$ and $(2)$), and hence, $u \equiv 0$ on $M$ by the unique continuation theorem; for details and other interesting facts, see \cite{K5}. 
\end{proof}

\begin{lem}[{\bf precompactness}]
Let $\{ z_k \}_{k=1}^{\infty} \subset K_{\pm}$ be a sequence, and $\{ f_k \}_{k=1}^{\infty}$ be a bounded sequence in $L^2_{\frac{1}{2}+s}(M,\mu)$. 
Assume that $\{ u_k \}_{k=1}^{\infty}$ is the sequence of the corresponding outgoing or incoming solution to the equation $(10)$. 
If $\{ u_k \}_{k=1}^{\infty}$ is bounded in $L^2_{-\frac{1}{2}-s'}(M, \mu)$, then $\{ u_k \}_{k=1}^{\infty}$ is precompact in $L^2_{-\frac{1}{2}-s'}(M, \mu)$. 
\end{lem}
\begin{proof}
First, recall that $K_{+}$ and $K_{-}$ are any fixed compact subsets in ${\it \Pi}_{+}$ and  ${\it \Pi}_{-}$, respectively. 
Therefore, 
\begin{align}
  \min \Big\{ \,\Big|z - \frac{(\beta_j)^2}{4} \Big|~\Big|~ z \in K_{\pm},~j=1,\cdots,m \Big\} >0. 
\end{align}
For $\varepsilon >0$, let ${\rm Rect}_{\pm} \big( \frac{(\beta_j)^2}{4}, \varepsilon \big)$ be a rectangle around $\frac{(\beta_j)^2}{4}$ defined by 
\begin{align*}
 {\rm Rect}_{\pm} \Big( \frac{(\beta_j)^2}{4}, \varepsilon \Big) 
 := \Big\{ x + iy \in {\it \Pi}_{\pm} \,\big|\, 
 \max \Big\{ |y|,\,\Big|\frac{(\beta_j)^2}{4} - x \Big| \Big\} \le \varepsilon \Big\}.
\end{align*}
Then, $(71)$ implies that there exists a constant $\varepsilon_0 = \varepsilon_0 (K_{\pm}, c_1, \cdots, c_m) > 0$ such that, ${\rm Rect}_{+} \big( \frac{(\beta_j)^2}{4}, \varepsilon_0 \big) \cap K_{+} = \emptyset$ for any $j=1, \cdots, m$, and ${\rm Rect}_{-}\big( \frac{(\beta_j)^2}{4}, \varepsilon_0 \big) \cap K_{-} = \emptyset$ for any $j=1, \cdots, m$. 
Thus, for each fixed $j \in \{ 1, \cdots, m \}$, if $z\in K_{\pm}$, then one of the following three cases occurs: (a) ${\rm Re}\,z \le \frac{(\beta_j)^2}{4} - \varepsilon_0$; (b) $|{\rm Im}\, z| \ge \varepsilon_0$; (c) ${\rm Re}\,z \ge \frac{(\beta_j)^2}{4} + \varepsilon_0$. 
For the case (a), we shall apply Proposition $6.1$; 
for the case (b), we shall apply Corollary $8.2$; 
for the case (c), we shall apply Proposition $7.1$. 
Then, for any $R > \max\{1, r_0\}$, we obtain 
\begin{align}
 \|u_k\|_{L^2_{-\frac{1}{2}-s'}(E_j(R,\infty),\mu)} 
 \le 
 \widehat{c}_{12} \cdot (1+R)^{-s'} 
 \left\{ \| f_k \|_{L^2_{\frac{1}{2} + s}(M, \mu)} 
 + \| u_k \|_{L^2_{-\frac{1}{2}-s'}(M, \mu)} \right\},
\end{align}
where $\widehat{c}_{12}>0$ is a constant independent of $\{ u_k \}_{i=1}^{\infty}$ and $\{ f_k \}_{i=1}^{\infty}$. 
Since $(72)$ holds for any $j \in \{ 1, \cdots, m \}$, by summing up for every $j \in \{ 1, \cdots, m \}$, we obtain, for any $R > \max \{ 1, r_0 \}$, 
\begin{align}
  \|u_k\|_{L^2_{-\frac{1}{2}-s'}(E(R,\infty),\mu)} 
 \le 
  \widehat{c}_{14} \cdot (1+R)^{-s'} 
  \left\{ \| f_k \|_{L^2_{\frac{1}{2}+s}(M,\mu)} 
  + \| u_k \|_{L^2_{-\frac{1}{2}-s'}(M,\mu)} \right\},
\end{align}
where $\widehat{c}_{14}>0$ is a constant independent of $\{ u_k \}_{i=1}^{\infty}$ and $\{ f_k \}_{i=1}^{\infty}$. 
Note that the right hand side of $(73)$ tends to zero as $R\to \infty$ uniformly with respect to $k$. 

On the other hand, since $- L u_k = z_k u_k + f_k$, we see that $\{-Lu_k\}_{k=1}^{\infty}$ is a bounded sequence in $L^2_{-\frac{1}{2}-s'}(M,\mu)$, and hence, by Corollary $4.1$, $\bigl\{|\nabla u_k|\bigr\}_{k=1}^{\infty}$ is also a bounded sequence in $L^2_{-\frac{1}{2}-s'}(M,\mu)$. 
Therefore, the Rellich's lemma implies that 
\begin{align}
  \big\{ u_k|_{U(R)} \big\}_{k=1}^{\infty}~{\rm is~precompact~in}~
  L^2_{-\frac{1}{2}-s'}(U(R),\mu), \quad {\rm for~any}~R>0.
\end{align}
Thus, by $(73)$ and $(74)$, we see that $\{ u_k \}_{k=1}^{\infty}$ is precompact in $L^2_{-\frac{1}{2}-s'}(M,\mu)$. 
\end{proof}
\begin{lem}[{\bf preservation of radiation condition}]
Let $\{ z_k \}_{k=1}^{\infty} \subset K_{\pm}$ and $\{ f_k \}_{k=1}^{\infty} \subset L^2_{\frac{1}{2}+s}(M,\mu)$ be sequences. 
Assume that 
\begin{align*}
  \lim_{k\to \infty} z_k = z_{\infty}; \quad 
  f_k \to f_{\infty}~(k \to \infty) \quad 
  {\rm weakly~in~} L^2_{\frac{1}{2}+s}(M,\mu)
\end{align*}
and that $u_k$ is the outgoing or incoming solution of the equation $(10)$ with $z = z_k$ and $f=f_k$. 
Assume that $\displaystyle \lim_{k\to \infty}\|u_k-u_{\infty}\|_{L^2_{-\frac{1}{2}-s'}(M,\mu)}=0$. 
Then, $u_{\infty}$ is the outgoing or incoming solution of $(10)$ with $z=z_{\infty}$ and $f=f_{\infty}$. 
\end{lem}
\begin{proof}
By taking the limit of the equation 
\begin{align}
  - L u_k - z_k u_k = f_k,
\end{align}
we have 
\begin{align}
  - L u_{\infty} - z_{\infty} u_{\infty} = f_{\infty}
\end{align}
in the sense of distribution. 
But, the elliptic regularity theorem implies that $u_{\infty}\in H^2_{loc}(M)$, and hence, the equation $(76)$ holds in the sense of $L^2_{loc}(M)$. 
On the other hand, since any weakly convergent sequence is bounded by the principle of uniform boundedness, $\| f_k \|_{L^2_{\frac{1}{2}+s}(M,\mu)}$ is uniformly bounded with respect to $k$. 
Therefore, Proposition $5.2$ implies that $\| \nabla u_k + u_k \nabla p_{\pm}(z_k, *) \|^2_{L^2_{-\frac{1}{2}+s}(E(3,\infty),\mu)}$ is also uniformly bounded with respect to $k$. 
Hence, by taking a subsequence if necessary, we may assume that 
\begin{align}
  \nabla u_k + u_k \nabla p_{\pm} (z_k, *) \to X~(k\to \infty) \quad 
  {\rm weakly~in}~L^2_{-\frac{1}{2}+s}(E(3,\infty),\mu).
\end{align}
Since $\| f_k \|_{L^2_{\frac{1}{2}+s}(M,\mu)}$ and $\| u_k \|_{L^2_{-\frac{1}{2}-s'}(M,\mu)}$ are bounded uniformly with respect to $k$, the equation $(75)$ implies that, for any $R>0$, there exists a constant $\widehat{c}_0(R)$ such that $\|Lu_k\|_{L^2(U(R),\mu)} + \| u_k \|_{L^2(U(R),\mu)} \le \widehat{c}_0(R)$, where $\widehat{c}_0(R)$ is independent of $k$. 
Therefore $\| u_k \|_{H^2(U(R))}$ is uniformly bounded with respect to $k$. 
Hence, the Rellich's lemma implies that $\{ u_k|_{U(R)} \}_{k=1}^{\infty}$ is precompact in $H^1(U(R))$, and hence, our assumption, $\displaystyle\lim_{k\to \infty}\| u_k - u_{\infty} \|_{L^2_{-\frac{1}{2}-s'}(M,\mu)}=0$, implies that
\begin{align}
  \lim_{k\to \infty} \| u_k - u_{\infty} \|_{W^{1,1}(U(R),\mu)} = 0 \quad 
  {\rm for~any}~R > 0.
\end{align}
In view of $(77)$ and $(78)$, we see that $X \in L^2_{-\frac{1}{2}+s}(E(3,\infty),\mu)$ coincides with $\nabla u_{\infty}+u_{\infty}\nabla p_{\pm}(z_{\infty},\cdot)$, and hence, $u_{\infty}$ satisfies the condition $(14)$. 
Thus, we have proved Lemma $8.5$. 
\end{proof}
\begin{lem}[{\bf boundedness of solutions}]
For $(z,f)\in K_{\pm}\times L^2_{\frac{1}{2}+s}(M, \mu)$, let $u(z,f)$ denote the corresponding outgoing or incoming solution $u(z,f)$ of $(10)$. 
Then, there exists a constant $c(K_{\pm})>0$, depending only on $K_{\pm}$, such that 
\begin{align*}
  \|u(z, f)\|_{L^2_{-\frac{1}{2}-s'}(M,\mu)} 
\le 
  c(K_{\pm}) \, \|f\|_{L^2_{\frac{1}{2}+s}(M,\mu)}.
\end{align*}
\end{lem}
\begin{proof}
We shall prove this lemma by contradiction. 
If we deny the conclusion, there exist sequences $\{z_k\}_{k=1}^{\infty} \subset K_{\pm}$ and $\{f_k\}_{k=1}^{\infty} \subset L^2_{\frac{1}{2}+s}(M,\mu)$ such that 
\begin{align*}
  \|u_k\|_{L^2_{-\frac{1}{2}-s'}(M,\mu)} \equiv 1; \quad 
  \lim_{k \to \infty} \|f_k\|_{L^2_{\frac{1}{2}+s}(M,\mu)} =0; \quad 
  \lim_{k \to \infty} z_k = z_{\infty} \in K{\pm},
\end{align*}
where $u_k:=u(z_k, f_k)$ and we have used the fact that $K_{\pm}$ is compact. 
By Lemma $8.2$, by taking a subsequence if necessary, we may assume that there exists $u_{\infty} \in L^2_{-\frac{1}{2}-s'}(M,\mu)$ such that $\displaystyle\lim_{k\to \infty} \| u_k-u_{\infty} \|_{L^2_{-\frac{1}{2}-s'}(M,\mu)} = 0$. 
Then, by taking the limit of the equation $- L u_k - z_k u_k = f_k$, we obtain 
\begin{align}
  - L u_{\infty} - z_{\infty} u_{\infty} = 0 
\end{align}
in the sense of distribution.
But, the elliptic regularity theorem implies that $u\in C^{\infty}(M)$ and $(79)$ holds in the sense of $C^{\infty}(M)$. 
Hence, if ${\rm Im}\,z_{\infty} \neq 0$, Corollary $8.1$ and Proposition $8.1$ (i) imply that $u_{\infty} \equiv 0$; if ${\rm Im}\,z_{\infty} = 0$, Lemma $8.3$ and Lemma $8.1$ imply that $u_{\infty} \equiv 0$. 
This contradicts the fact that $\|u_{\infty}\|_{L^2_{-\frac{1}{2}-s'}(M,\mu)} 
= \lim_{k\to \infty} \|u_k\|_{L^2_{-\frac{1}{2}-s'}(M,\mu)} = 1$. 
Thus, we have proved Lemma $8.4$ 
\end{proof}
We are now ready to prove the following:
\begin{thm}[{\bf principle of limiting absorption}]
Let $s$ and $s'$ be constants satisfying $(59)$. 
Then, in the Banach space $\mathbb{B}_{\mu}(\frac{1}{2}+s,-\frac{1}{2}-s')$, we have the limit 
\begin{align}
  R_{-L}(\lambda \pm i\,0)
  = \lim_{\varepsilon \downarrow 0} R_{-L} (\lambda \pm i \varepsilon),
  \qquad \lambda \in I.
\end{align}
Moreover, this convergence $(80)$ is uniform on any compact subset of $I$, and $R_{-L}(z)$ is continuous on ${\it \Pi}_{+}$ and ${\it \Pi}_{-}$ with respect to the operator norm $\| * \|_{\mathbb{B}_{\mu}(\frac{1}{2}+s,-\frac{1}{2}-s')}$ by considering $R_{-L}(\lambda +i\,0)$ and $R_{-L}(\lambda -i\,0)$ on ${\it \Pi}_{+}\cap (0,\infty)$ and ${\it \Pi}_{-}\cap (0,\infty)$, respectively. 
\end{thm}
\begin{proof}
Let $K$ be any compact subset of $I$, and we will show that there exists an operator $R_{-L}(\lambda+i0)$ in $\mathbb{B}_{\mu}(\frac{1}{2}+s,-\frac{1}{2}-s')$ such that 
\begin{equation}
 \lim_{\tau \downarrow 0} \sup_{\lambda\in K} 
 \| R_{-L}(\lambda + i\,0) - R_{-L}(\lambda + i\tau) \|
 _{\mathbb{B}_{\mu}(\frac{1}{2}+s,-\frac{1}{2}-s')} = 0
\end{equation}
by contradiction. 
If we assume the contrary, there exist constant $\varepsilon_0 > 0$, sequences $\{ \tau_k^1 \}_{k=1}^{\infty}$, $\{ \tau_k^2 \}_{k=1}^{\infty}$ of positive numbers, $\{f_k \}_{k=1}^{\infty} \subset L^2_{\frac{1}{2}+s}(M,\mu)$, and $\{\lambda_k\}_{k=1}^{\infty} \subset K$ such that 
\begin{align}
 & 0 < \tau_k^1, \tau_k^2 < 1/k ~;~
   \|f_k\|_{L^2_{\frac{1}{2}+s}(M,\mu)} = 1~; \notag \\
 & \| u(\lambda_k + i\tau_k^1,f_k) - u(\lambda_k + i\tau_k^2,f_k) \|
   _{L^2_{-\frac{1}{2}-s'}(M,\mu)} \ge \varepsilon_0,
\end{align}
where we have used the notation in Lemma $8.4$. 
By taking a subsequence if necessary, we may assume that 
\begin{align*}
  \lim_{k\to \infty}\lambda_k = \lambda_{\infty}~;~
  f_k\to f_{\infty}~(k\to \infty) 
  \quad {\rm weakly~in}~L^2_{\frac{1}{2}+s}(M,\mu)~; 
\end{align*}
moreover, by Lemma $8.4$ and Lemma $8.2$, 
\begin{align}
& \lim_{k\to \infty} \| u(\lambda_k + i\tau_k^1, f_k) - u_{\infty} \|
   _{L^2_{-\frac{1}{2}-s'}(M,\mu)} =0~; \\
&  \lim_{k\to \infty} \| u(\lambda_k + i\tau_k^2, f_k) - u_{\infty}' \|_{L^2_{-\frac{1}{2}-s'}(M,\mu)} = 0. 
\end{align}
Then, Lemma $8.3$ implies that $u_{\infty}$ and $u_{\infty}'$ are outgoing solutions of the same equation $(10)$ with $z = \lambda_{\infty}$ and $f = f_{\infty}$. 
However, Lemma $8.1$ implies that $u_{\infty}=u_{\infty}'$, which contradicts $(82)$, $(83)$, and $(84)$. 
Thus we have proved $(81)$. 

The proof of the existence of an operator $R_{-L}(\lambda - i\,0) \in \mathbb{B}_{\mu}(\frac{1}{2}+s, -\frac{1}{2}-s')$ satisfying $\lim_{\tau \downarrow 0} \sup_{\lambda\in K} \| R_{-L}(\lambda - i\,0) - R_{-L}(\lambda-i\tau)\|_{\mathbb{B}_{\mu}(\frac{1}{2}+s,-\frac{1}{2}-s')} = 0$ is quite the same. 
\end{proof}
\begin{cor}[{\bf absolutely continuity}]
Let $(M,g)$ be a Riemannian manifold as in Theorem $1.1$. 
Then, $-L$ is absolutely continuous on $\big( \frac{(\beta_1)^2}{4}, \infty \big)$ and has no singular continuous spectrum. 
\end{cor}
\begin{proof}
First, note that the essential spectrum of $-L$ is equal to $\big[ \frac{(\beta_1)^2}{4}, \infty \big)$ (see \cite{K1}, \cite{K2}), because $\Delta_g r \to \beta_j$ as $r \to \infty$ on $E_j$ for $j=1,2, \cdots, m$. 

It is now a standard fact that the limiting absorption principle implies the absolute continuity; see [28, Theorem XIII.$19$ and Theorem XIII.$20$]. 
Hence, we obtain the absolute continuity of $-L$ on $I = \big(\frac{(\beta_1)^2}{4},\infty \big) - \left\{ \frac{(\beta_j)^2}{4} \mid j=2, \cdots, m \right\}$ by Theorem $8.1$. 
Moreover, Lemma $8.1$ implies that $\frac{(\beta_j)^2}{4} \notin \sigma_{\rm pp}(-L)$ for $j=2, \cdots, m$ satisfying $\beta_j > \beta_1$. 
Here, $\sigma_{\rm pp}(-L)$ stands for the point spectrum of $-L$. 
Hence, $-L$ is absolutely continuous on $\big( \frac{(\beta_1)^2}{4}, \infty \big)$, and has no singular continuous spectrum. 
\end{proof}
\begin{cor}
Let $(M,g)$ be a Riemannian manifold as in Theorem $8.1$. 
When $M$ has at least one end satisfying {\rm MC}$(\frac{a}{r}, \frac{b}{r}, \delta)$'' for some $a>0$, $b>0$, and $\delta>0$, then $-L$ is absolutely continuous on $(0, \infty)$ and $0\notin \sigma_{\rm pp}(-L)$. 
\end{cor}
Since the multiplication operator $e^w : L^2_{\alpha}(M,v_M) \to L^2_{\alpha}(M,\mu)$ is unitary for any $\alpha \in \mathbb{R}$, operators, $-L$ and $-\Delta_g$, are unitarily equivalent. 
Hence, Theorem $8.1$, Corollary $8.3$, and Corollary $8.4$ imply Theorem $1.1$, Theorem $1.2$, and Corollary $1.1$, respectively.

\section{Further discussions}

Corollary $9.1$ below means that ends into which a ``wave function'' $e^{it\Delta_g}u$ will recede as $t\to \pm \infty$ are not high energy ends $\cup \left\{ E_k \mid \beta_k \ge \beta_j,~k = 2, \cdots, m \right\}$ but low energy ends $\cup \left\{ E_k \mid \beta_k < \beta_j,~k=1, \cdots, j-1 \right\}$, when $u \in E_{-\Delta_g}(I_j) L^2(M, v_g)$ and $I_j = \big( \frac{(\beta_1)^2}{4}, \frac{(\beta_j)^2}{4} \big)$. 

\begin{cor}
Let $(M,g)$ be an $n$-dimensional Riemannian manifold as in Theorem $1.1$, and 
$E_{-\Delta_g}(\Lambda)$ $(\Lambda \in \mathcal{B})$ denotes the spectral resolution of $-\Delta_g$ on $L^2(M,v_g)$. 
Let $j \in \{2,\cdots ,m\}$ be an integer satisfying $\beta_1 < \beta_j$, and set $I_j := \big( \frac{(\beta_1)^2}{4}, \frac{(\beta_j)^2}{4} \big)$. 
Then, for any $u \in E_{-\Delta_g}(I_j) L^2(M, v_g)$, we obtain 
\begin{align*}
  \lim_{t\to \pm \infty} \int_{\cup \left\{ E_k \mid \beta_k \ge \beta_j,~k=2,\cdots,m \right\} } |e^{it\Delta}u|^2 \,dv_g = 0.
\end{align*}
\end{cor}
\begin{proof}
For any constant $\varepsilon >0$ satisfying $\frac{(\beta_1)^2}{4} < \frac{(\beta_j)^2}{4} - \varepsilon$, we donote $I_j(\varepsilon) := \big( \frac{(\beta_1)^2}{4}, \frac{(\beta_j)^2}{4} - \varepsilon \big)$. 
Then, $E_{-\Delta_g} \left( I_j(\varepsilon) \right)$ strongly converges to $E_{-\Delta_g} \left( I_j \right)$ as $\varepsilon \to +0$. 
Hence, it suffices to consider the case that $u = E_{-\Delta_g} \left( I_j(\varepsilon) \right) u$. 
Theorem $1.2$ implies that $u$ is an element in the absolutely continuous subspace of $-\Delta_g$; hence, $(E_{-\Delta_g}(\Lambda)u, v)$ $(\Lambda \in \mathcal{B})$ is an absolutely continuous signed measure, for any $v \in L^2(M, v_g)$. 
Thus, there exists a function $f(\lambda) \in L^1(\mathbb{R})$ such that $( E_{-\Delta_g}(\Lambda) u,v)_{L^2(M, v_g)} = \int_{\Lambda} f(\lambda) \,d\lambda$ for any $\Lambda \in \mathcal{B}$. 
Therefore, for each integer $k \ge 0$, we obtain
\begin{align}
& \big( e^{it\Delta_g} (\Delta_g)^k (1-\Delta_g) u, v \big)_{L^2(M, v_g)} 
  = \int_{I_j(\varepsilon)} e^{it\lambda} \lambda^k (1+\lambda) 
  \, d(E_{-\Delta_g}(\lambda)u, v) \\
= 
& \int_{I_j(\varepsilon)} e^{it\lambda} \lambda^k (1 + \lambda) f(\lambda) \,d\lambda \notag~;
\end{align}
Riemann--Lebesgue lemma implies that the last term of $(85)$ converge to zero as $t\to \pm$. 
Thus, $e^{it\Delta_g} (\Delta_g)^k (1-\Delta_g) u$ weakly converges to zero as $t\to \pm\infty$ in $L^2(M, v_g)$ for each $k \ge 0$. 
Set
\begin{align*}
 \chi_{\displaystyle_{U(R)}} (x) := 
 \begin{cases}
  \ 1 \quad & {\rm if}~~x \in U(R),\\
  \ 0 \quad & {\rm if}~~x \in M \backslash U(R).
 \end{cases}
\end{align*}
Then, the Rellich's lemma implies that $\chi_{\displaystyle_{U(R)}} (1 - \Delta_g)^{-1}$ is a compact operator on $L^2(M, v_g)$. 
Therefore, $\chi_{\displaystyle_{U(R)}} (\Delta_g)^k e^{it\Delta_g}u = \chi_{\displaystyle_{U(R)}} (1-\Delta_g)^{-1} e^{it\Delta_g} (\Delta_g)^k (1 - \Delta_g) u$ strongly converges to zero as $t\to \pm \infty$ in $L^2(M, v_g)$, that is, 
$$
  \lim_{t\to \pm \infty} \int_{U(R)} |(\Delta_g)^k e^{it\Delta_g}u|^2 \,dv_g 
  = 0
$$
for any nonnegative integer $k$ and $R>0$. 
Thus, the Sobolev embedding theorem implies that, for any integer $k \ge 0$ and  $R>0$, 
\begin{align}
  e^{it\Delta_g} u \to 0 \quad {\rm as}~t\to \pm \infty \quad {\rm on}~U(R) 
\quad {\rm in~the~sense~of}~C^k\mbox{--topology}.
\end{align}
Now, we shall prove
\begin{align}
  \lim_{t\to \infty} 
  \int_{ \cup \left\{ E_k \mid \beta_k \ge \beta_j,~k =2, \cdots, m \right\} } 
  |e^{it\Delta_g}u|^2 \,dv_g = 0
\end{align}
by contradiction. 
If we assume contrary, there exist $j_1 \in \{ 2, \cdots, m \}$ satisfying $\beta_{j_1} \ge \beta_j$, a positive constant $a_{j_1}$, and a divergent sequence $\{t_{\ell}\}_{\ell =1}^{\infty}$ of positive real numbers such that 
\begin{align}
  \lim_{\ell\to \infty} \int_{E_{j_1}} |e^{it_{\ell}\Delta_g}u|^2 \,dv_g 
  \ge a_{j_1} > 0.
\end{align}
Then, since
\begin{align*}
  \inf \bigg\{ 
  \frac{\| \nabla h \|_{L^2(E_{j_1}(R,\infty), v_g)}^2}{\| h \|_{L^2(E_{j_1}(R,\infty), v_g)}^2} ~\bigg|~0 \neq h \in C_0^{\infty}\big( E_{j_1}(R,\infty) \big) \bigg\} \to \frac{( \beta_{j_1} )^2}{4} \quad {\rm as}~R\to \infty,
\end{align*}
we see that $(86)$ and $(88)$ yield
\begin{align}
& \lim_{\ell \to \infty} 
  \| (\Delta_g)^k e^{it_{\ell}\Delta_g}u \|_{L^2(M, v_g)}^2
=
  \sum_{j=1}^{m} \lim_{\ell \to \infty} 
  \| (\Delta_g)^k e^{it_{\ell}\Delta_g}u \|_{L^2(E_j, v_g)}^2 \\
\ge 
 & \lim_{\ell \to \infty} \| (\Delta_g)^k e^{it_{\ell}\Delta_g}u \|
   _{L^2(E_{j_1},v_g)}^2 \ge \Big( \frac{\beta_{j_1}}{2} \Big)^{4k} a_{j_1} . \notag 
\end{align}
On the other hand, since $u = E \left( I_j(\varepsilon) \right) u$, 
we have, for any $t\in \mathbb{R}$, 
\begin{align}
& \| (\Delta_g)^ke^{it\Delta_g}u \|_{L^2(M, v_g)}^2
= \int_{I_j(\varepsilon)} \lambda^{2k} \,d\| E_{-\Delta_g}(\lambda)u\|^2 \\
\le 
 & \Big( \frac{(\beta_j)^2}{4} - \varepsilon \Big)^{2k} 
   \int_{I_j(\varepsilon)} \,d \| E_{-\Delta_g}(\lambda)u \|^2 
= \Big( \frac{(\beta_j)^2}{4} - \varepsilon \Big)^{2k} \| u \|_{L^2(M,v_g)}^2. \notag 
\end{align}
From $(89)$ and $(90)$, we obtain $\theta^{2k} \| u \|_{L^2(M,v_g)}^2 \ge a_{j_1}$ for all integer $k \ge 0$, where $\theta := \big( (\beta_j)^2 - 4 \varepsilon \big)/(\beta_j)^2$; since $\theta \in (0,1)$, letting $k\to \infty$, we obtain $a_{j_1}=0$, which contradicts $(88)$. 
This completes the proof of $(87)$. 
The proof of the case, $t\to -\infty$, is quite the same. 
\end{proof}

In view of Proposition $6.1$, we see that the following holds:
\begin{thm}
Let $(M,g)$ be an $n$-dimensional connected complete noncompact Riemannian manifold and $U$ be a relatively compact open subset of $(M, g)$. 
Assume that $M \backslash U$ consists of the disjoint union of ends, $E_1, \cdots, E_m$, with radial coordinates, where $m \ge 2$. 
Assume that there exist positive constants, $\gamma $ and $\delta$, such that
\begin{align*}
 & E_j~{\rm satisfies~MC} \Big( \frac{a_j}{r}, \frac{b_j}{r}, \delta \Big)
   \quad {\rm for}~ 1 \le j \le m_0~;\\
 & E_j~{\rm satisfies~MC}(\alpha_j, \beta_j, \delta)
   \quad \ \,\,{\rm for}~~m_0 + 1 \le j \le m_1~;\\
 & \lim_{t\to \infty} \inf \big\{ \Delta_g r(x) \mid x \in E_{m_1+1}(t,\infty) 
   \cup \cdots \cup E_m(t,\infty) \big\} \ge \gamma > \frac{(\beta_1)^2}{4},
\end{align*}
where 
\begin{align*}
  0 = \beta_1 = \beta_2 = \cdots = \beta_{m_0} < \beta_{m_0+1} 
  \le \beta_{m_0+2} \le \cdots \le \beta_{m_1}
\end{align*}
are real constants\,{\rm ;} $0 \le m_0 \le m_1$ and $1 \le m_1 \le m-1$ are integer\,{\rm ;} $a_j,b_j,\alpha_j,\beta_j$ are all positive constants. 
Then, in the Banach space $\mathbb{B}(\frac{1}{2}+s,-\frac{1}{2}-s')$, 
we have the limit 
\begin{equation*}
  R(\lambda \pm i0)
:= \lim_{\varepsilon \downarrow 0}R(\lambda \pm i\varepsilon), \qquad 
  \lambda \in I':= \Big( \frac{(\beta_1)^2}{4}, \gamma \Big) 
  - \Big\{ \frac{(\beta_k)^2}{4} \, \big| \, 1 \le k \le m_1 \Big\},
\end{equation*}
where $s$ and $s'$ are constants satisfying $(59)$. 
Moreover, this convergence is uniform on any compact subset of $I'$ and $-\Delta_g$ is absolutely continuous on $\big( \frac{(\beta_1)^2}{4}, \gamma \big)$. 
\end{thm}
\begin{proof}
Reconsidering the arguments in Section $8$, we see that the precompactness lemma plays an important role there; it follows from decay estimates of solutions of $(10)$ on every ends. 
Hence, when $\lambda \in I'$, we shall apply Proposition $6.1$ for ends $E_{m_1+1}, \cdots, E_m$; for the rest of ends, we shall apply the arguments in Section $8$. 
Then, we obtain the precompactness lemma. 
The rest of the proof holds good, and we obtain Theorem $9.1$. 
\end{proof}


\vspace{10mm}

\end{document}